\documentclass[twoside]{amsart}
\usepackage{latexsym}
\usepackage{amssymb,amsmath,amsopn}
\usepackage[dvips]{graphicx}   %To insert ps figures
\usepackage{color,epsfig}      %To insert ps figures 

\newtheorem{thm}{Theorem}[section]
\newtheorem{prop}[thm]{Proposition}
\newtheorem{claim}[thm]{Claim}
\newtheorem{cor}[thm]{Corollary}
\newtheorem{lem}[thm]{Lemma}
\newtheorem{conj}[thm]{Conjecture}
\newtheorem{obs}[thm]{Observation}
\newtheorem{prob}[thm]{Problem}
\theoremstyle{definition}
\newtheorem{defn}[thm]{Definition} 
\newtheorem{rem}[thm]{Remark}
\newtheorem{ex}[thm]{Example}
\numberwithin{equation}{section}

\renewcommand{\Re}{\mathbb R}

\renewcommand{\S}{\mathbb S}
\newcommand{\Ze}{\mathbb Z}

\newcommand{\B}{\mathbf{B}}

\newcommand{\fac}{\mathcal{F}}

\newcommand{\vv}{v}
\newcommand{\vu}{u}
\newcommand{\Sph}{\S^{n-1}(o)}
\newcommand{\F}{\mathfrak{F}}

\DeclareMathOperator{\scsc}{Sconv}
\newcommand{\Sconv}[2]{\scsc\left(#1,#2\right)}
\DeclareMathOperator{\inter}{int}
\DeclareMathOperator{\bd}{bd}
\DeclareMathOperator{\diam}{diam}
\DeclareMathOperator{\conv}{conv}
\DeclareMathOperator{\aff}{aff}
\DeclareMathOperator{\card}{card}

\DeclareMathOperator{\area}{Area}
\DeclareMathOperator{\convv}{\conv_s}
\DeclareMathOperator{\per}{Perimeter}
\DeclareMathOperator{\crr}{cr}

\begin{document}
\title[ball-polyhedra]{ball-polyhedra}
\author[K. Bezdek, Z. L\'angi, M. Nasz\'odi \and P. Papez]{by\\\;\\
K\'aroly Bezdek$^{*}$, Zsolt L\'angi$^{**}$,\\\;\\
M\'arton Nasz\'odi$^{**}$ and Peter Papez$^{***}$}

\address{K\'aroly Bezdek, Dept. of Math. and Stats., University of Calgary,
2500 University Drive NW, Calgary, Ab, Canada T2N 1N4}
\email{bezdek@math.ucalgary.ca}
\address{Zsolt L\'angi, Dept. of Math. and Stats., University of Calgary,
2500 University Drive NW, Calgary, Ab, Canada T2N 1N4}
\email{zlangi@math.ucalgary.ca}
\address{M\'arton Nasz\'odi, Dept. of Math. and Stats., University of Calgary,
2500 University Drive NW, Calgary, Ab, Canada T2N 1N4}
\email{nmarton@math.ucalgary.ca}
\address{Peter Papez, Dept. of Math. and Stats., University of Calgary,
2500 University Drive NW, Calgary, Ab, Canada T2N 1N4}
\email{pdpapez@math.ucalgary.ca}
                                                                               
\thanks{$*$ Partially supported by the Hung. Nat. Sci. Found.
(OTKA), grant no. T043556 and T037752 and by a Natural Sciences and Engineering 
Research Council of Canada Discovery Grant.}
\thanks{$**$ Partially supported by the Hung. Nat. Sci. Found.
(OTKA), grant no. T043556 and T037752, and by the Alberta Ingenuity Fund.}
\thanks{$***$ Partially supported by an Alberta Graduate Fellowship.}


\begin{abstract}
We study two notions. One is that of
\emph{spindle convexity}. A set of circumradius not greater than one is
spindle convex if, for any pair of its points, it contains every short
circular arc of radius at least one, connecting them. The other objects of study are
bodies obtained as intersections of finitely many balls of the same
radius, called \emph{ball-polyhedra}. We find analogues of several
results on convex polyhedral sets for ball-polyhedra. 
\end{abstract}
\maketitle

\section{Introduction}

The main goal of this paper is to study the geometry of intersections of
finitely many  congruent balls, say of unit balls, from the viewpoint of
discrete geometry in Euclidean space.  We call these sets
\emph{ball-polyhedra}. Some special classes have been studied in the
past; see, e.g. \cite{Bieberbach1}, \cite{Bieberbach2} and \cite{Mayer}. 
For Reuleaux polygons see \cite{KupMar} and \cite{KMW}. Nevertheless, the
name ball-polyhedra  seems to be a new terminology for this special
class of linearly convex sets. In fact, there is a special kind of
convexity  entering along with ball-polyhedra which we call
\emph{spindle convexity}. We thank the referee for suggesting this name
for this notion of convexity that was first introduced by Mayer
\cite{Mayer} as ``\"Uberkonvexit\"at''.

The starting point of our research described in this paper was a sequence of
lectures of the first named author on ball-polyhedra given at the
University of Calgary in the fall of 2004. Those lectures have been
strongly motivated by  the following recent papers that proved
important new geometric properties of  intersections of finitely many
congruent balls: a proof of the Borsuk conjecture for finite point sets in three-space
based on the combinatorial geometry of ``spherical polytopes'' 
(\cite{AP}, p. 215); Sallee's theorem \cite{Salle}
claiming that the class of the so-called ``Reuleaux polytopes'' is dense
in the class of sets of constant width in $\Re^3$; a proof of the
Kneser-Poulsen conjecture in the plane by K. Bezdek and Connelly
\cite{BeCo} including the claim that under any contraction of the
center points of finitely many circular disks of $\Re^2$ the area of the
intersection cannot decrease, and finally an analogue of Cauchy's
rigidity theorem for triangulated ball-polyhedra in $\Re^3$ \cite{BN}. In
addition it should be noticed that ball-polyhedra play an essential role in 
the proof of Gr\"unbaum-Heppes-Straszewicz theorem on the maximal number of  
diameters of finite point sets in $\Re^3$; see \cite{KMP}.

This paper is not a survey on ball-polyhedra. Instead, it lays a rather
broad ground for future study of  ball-polyhedra by proving several new
properties of them and raising open research problems as well.

The structure of the paper is the following. First, notations
and basic results about spindle convex sets are introduced in
Sections~\ref{sec:notations} and \ref{sec:separation}. Some of these
results demonstrate the techniques that are different from the ones
applied in the classical theory. It seems natural that a more analytic
investigation of spindle convexity might belong to the realm of
differential geometry. 

In Section~\ref{sec:Kirchberger}, we find analogues of the theorem of
Kirchberger for separation by spheres. In
Section~\ref{sec:Caratheodory}, we prove spindle convex analogues of the
classical theorems of Carath\'eodory and Steinitz regarding the linear
convex hull.

In Section~\ref{sec:EPF}, we make the first steps in understanding the
boundary structures of ball-polyhedra. We present examples that show
that the face-structure of these objects is not at all obvious to
define. Section~\ref{sec:Maehara} contains our results on intersections
of unit spheres in $\Re^n$. The questions discussed there are motivated
primarily by a problem of Maehara \cite{Mae} and are related to the goal of
describing faces of ball-polyhedra. Also, we construct a counter-example to a
conjecture of Maehara in dimensions at least $4$. Then, in
Section~\ref{sec:monotonicity}, we discuss variants of the important
Kneser--Poulsen problem. In Section~\ref{sec:Steinitz}, we provide a partial characterization of
the edge-graphs of ball-polyhedra in $\Re^3$, similar to the Theorem of
Steinitz regarding convex polyhedra in $\Re^3$.

Then, in Section~\ref{sec:symmetricsections}, a conjecture of the first
named author about convex bodies in $\Re^3$ with axially symmetric sections is
proved for ball-polyhedra in $\Re^3$. We extend an illumination result
in $\Re^3$ of Lassak \cite{Lassak} and Weissbach \cite{Weissbach} in
Section~\ref{sec:illumination}. In Section~\ref{sec:Dowker}, we prove
various analogues of Dowker--type isoperimetric inequalities for
two-dimensional ball-polyhedra based on methods of L. Fejes-T\'oth
\cite{FT}. Finally, in Section~\ref{sec:Erdos-Szekeres}, we examine
spindle convex variants of Erd\H os--Szekeres-type questions.

\section{Notations and Some Basic Facts}\label{sec:notations}

Let $(\Re^n, ||\quad||)$, where $n \geq 2$, be the standard Euclidean
space with the usual  norm and denote the origin by $o$.  The Euclidean
distance between $a\in\Re^n$ and $b\in\Re^n$ is $||a-b||$. The closed 
line segment between two points is denoted by $[a,b]$, the open line
segment is  denoted by $(a,b)$. For the closed, $n$-dimensional ball
with center $a\in\Re^n$ and of radius $r>0$ we use the notation
$\B^n[a,r]:=\{x\in\Re^n:||a-x|| \leq r\}$. For the open $n$-dimensional
ball with center $a\in\Re^n$ and of radius $r>0$ we use the notation
$\B^n(a,r):=\{x\in\Re^n:||a-x||<r\}$. The $(n-1)$-dimensional  sphere
with center $a\in\Re^n$ and of radius $r>0$ is denoted by
$\S^{n-1}(a,r):=\{x\in\Re^n:||a-x||=r\}$. Any sphere or ball in the
paper is of positive radius. When $r$ is omitted, it is assumed to be
one. Using the usual conventions, let $\card,\conv,\inter,\bd$ and $\diam$
denote cardinality, convex hull, interior, boundary and
diameter of a set, respectively. We note that a $0$-dimensional
sphere is a pair of distinct points.

We introduce the following additional notations. For a set
$X\subset\Re^n$ let 

\begin{equation}
\B[X]:=\bigcap\limits_{x\in X}\B^n[x]\;\;\mbox{and}\;\;
\B(X):=\bigcap\limits_{x\in X}\B^n(x).
\end{equation}

\begin{defn}\label{defn:spindle}
Let $a$ and $b$ be two points in $\Re^n$. 
If $||a-b||<2$, then the \emph{closed spindle} of $a$ and $b$, 
denoted by $[a,b]_s$, is defined as the union of circular arcs  with endpoints 
$a$ and $b$ that are of radii at least one and are shorter than a
semicircle. If $||a-b||=2$, then $[a,b]_s:=\B^n[\frac{a+b}{2}]$.
If $||a-b||>2$, then we define $[a,b]_s$ to be $\Re^n$.

The \emph{open spindle}, denoted as $(a,b)_s$, in all cases is the interior of the closed one.
\end{defn}

\begin{rem}\label{rem:spindle}
If $||a-b||\leq 2$, then $[a,b]_s:=\B[\B[\{a,b\}]]$, and 
$(a,b)_s:=\B(\B[\{a,b\}])$.  
\end{rem}

\begin{defn}
The \emph{circumradius} $\crr(X)$ of a bounded set $X\subseteq\Re^n$ is defined as the radius of the unique smallest
ball that contains $X$ (also known as the circumball of $X$); that is,
\[
\crr(X):=\inf\{r>0: X\subseteq\B^n[q,r] \mbox{ for some } q\in\Re^n\}.
\] 
If $X$ is unbounded, then $\crr(X)=\infty$.
\end{defn}

Now, we are ready to introduce two basic notions that are used throughout this paper.

\begin{defn}
A set $C\subset\Re^n$ is \emph{spindle convex} if, for any pair of points
$a,b\in C$, we have $[a,b]_s\subseteq C$.
\end{defn}

\begin{defn}
Let $X\subset\Re^n$ be a finite set such that $\crr(X)\leq 1$.
Then we call $P:=\B[X]\neq\emptyset$ a \emph{ball-polyhedron}. For any $x\in X$ we call
$\B^n[x]$ a \emph{generating ball} of $P$ and $\S^{n-1}(x)$ a 
\emph{generating sphere} of $P$.
If $n=2$, then we call a ball-polyhedron a \emph{disk-polygon}.
\end{defn}

\begin{rem}
A spindle convex set is clearly convex. 
Moreover, since the spindle of two points has non-empty interior (if it exists),
a spindle convex set is either $0$-dimensional (if it is one point) or
full-dimensional. Also, the intersection of spindle convex sets is
again a spindle convex set.
\end{rem}

\begin{defn}
The \emph{arc-distance} of $a,b\in\Re^n$ is the 
arc-length of either shorter unit circular arcs connecting $a$ and $b$, when $||a-b|| \leq 2$; that is,
\[
\rho(a,b):=2\arcsin\left(\frac{\|a-b\|}{2}\right).
\]
If $||a-b|| > 2$, then $\rho(a,b)$ is undefined.
\end{defn}

\begin{rem}\label{rem:monotonicityofrho}
If $a,b,c\in\Re^n$ are points such that $||a-b||<||a-c||\leq 2$, 
then $\rho(a,b)<\rho(a,c).$
\end{rem}

The proof of the following claim is straightforward.
\begin{claim}[Euclidean arm-lemma]
Given two triangles with vertices $a,c,b$ and $a,c,b'$, respectively, in $\Re^2$ such
that $||c-b||=||c-b'||$ and the angle at $c$ in the first triangle is
less than in the second. Then $||a-b|| < ||a-b'||.$  
\end{claim}

In general, the arc-distance is not a metric. The following lemma
describes how the triangle-inequality holds or fails in some situations.
This lemma and the next corollary are from \cite{BeCoCs}, and they are
often applicable, as in Lemma~\ref{lem:concper}.

\begin{lem}
Let $a,b,c\in\Re^2$ be points such that 
$||a-b||, ||a-c||, ||b-c|| \leq 2.$ Then

\begin{tabular}{cccccl}
(i)  & $\rho(a,b)+\rho(b,c)$ & $>$ & 
$\rho(a,c)$ & $\Longleftrightarrow$ & $b\notin 
 [a,c]_s;$\\
(ii) &
$\rho(a,b)+\rho(b,c)$ & $=$ & $\rho(a,c)$ & $\Longleftrightarrow$ & $b\in

\bd[a,c]_s;$\\
(iii) & $\rho(a,b)+\rho(b,c)$ & $<$ & $\rho(a,c)$ &
$\Longleftrightarrow$ & $b\in  (a,c)_s.$
\end{tabular}
\end{lem}

\begin{cor}\label{cor:circlequadrilateral}
Let $a,b,c,d\in\Re^2$ be vertices of a spindle convex quadrilateral in
this cyclic order. Then
\[\rho(a,c)+\rho(b,d) > \rho(a,b)+\rho(c,d)\]
that is, the total arc-length of the diagonals is greater than the total
arc-length of an opposite pair of sides.
\end{cor} 

\section{Separation}\label{sec:separation}
This section describes results dealing with the separation of 
spindle convex sets by unit spheres motivated by the basic facts about 
separation of convex sets by hyperplanes as they are
introduced in standard textbooks; e.g., \cite{CR}.

\begin{lem}\label{lem:supportbyball}
Let a spindle convex set $C\subset\Re^n$ be supported by the 
hyperplane $H$ in $\Re^n$ at $x\in\bd C$. Then 
the closed unit ball supported by $H$ at $x$ and lying in the same side as $C$ contains $C$. 
\end{lem}

\begin{figure}[here]
\includegraphics[width=0.33\textwidth]{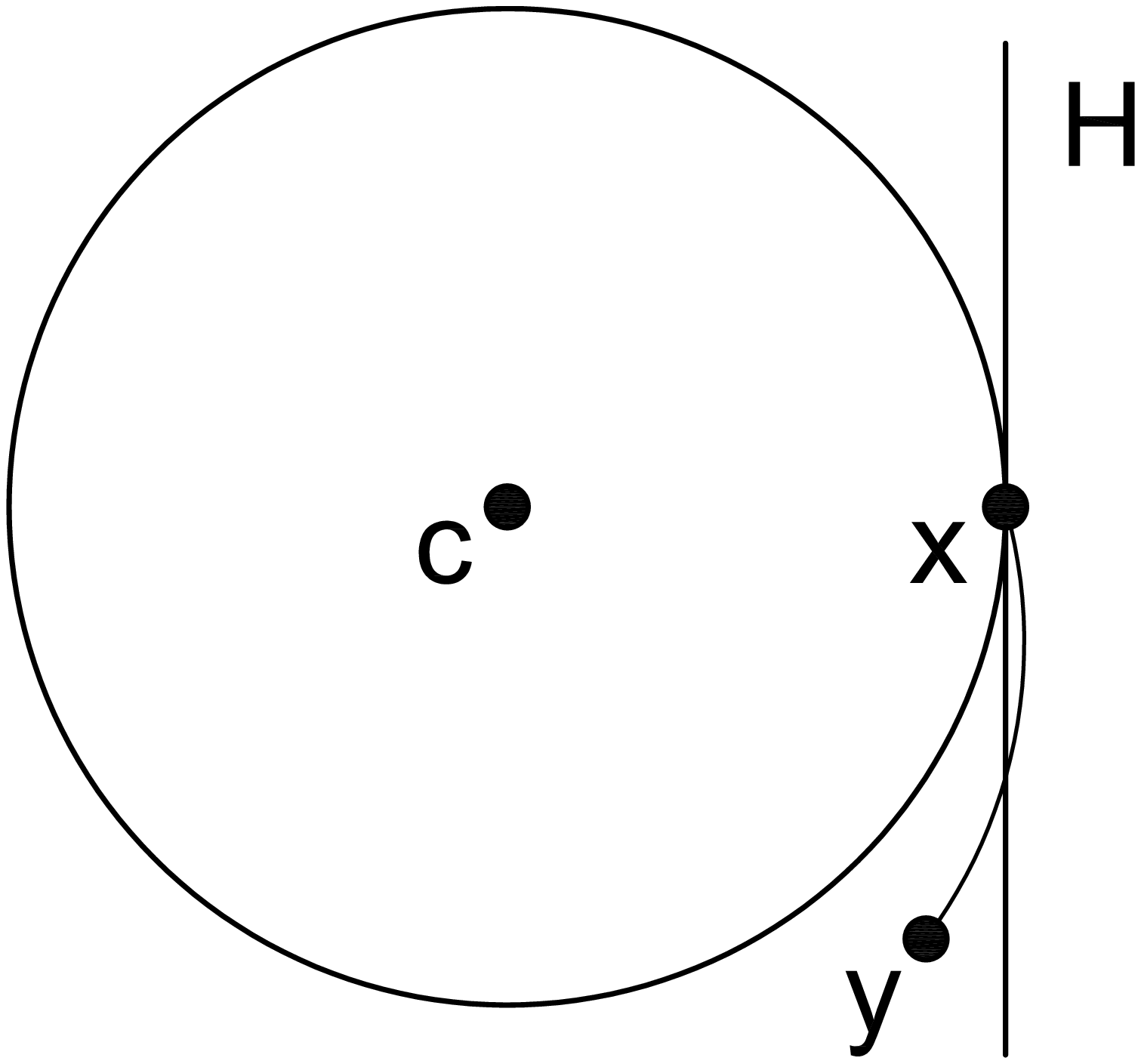}
\caption[]{}
\label{fig:separation}
\end{figure}

\begin{proof}
Let $\B^n[c]$ be the unit sphere that is supported by $H$ at $x$ 
and is in the same closed half-space bounded by $H$
as $C$. We show that $\B^n[c]$ is the desired unit ball. 

Assume that $C$ is not contained in $\B^n[c]$. So, there is a point
$y\in C,\;\;y\notin\B^n[c].$ Then, by taking the intersection
of the configuration with the plane that contains $x,y$ and $c$,
we see that there is a shorter unit circular arc connecting $x$ and
$y$ that does not intersect $\B^n(c)$ (Figure~\ref{fig:separation}). Hence, $H$ cannot be a
supporting hyperplane of $C$ at $x$, a contradiction.
\end{proof}

\begin{cor}\label{cor:crr}
Let $C\subset\Re^n$ be a spindle convex set.
If $\crr(C)=1$ then $C=\B^n[q]$ for some $q\in\Re^n$.
If $\crr(C)>1$ then $C=\Re^n$.  
\end{cor}

\begin{proof}
Observe that if $C$ has two distinct support unit balls then
$\crr(C)<1$. Thus, the first assertion follows. The second is clear.
\end{proof}

\begin{defn}
If a ball $\B^n[c]$ contains a set $C\subset\Re^n$ and a point
$x\in\bd C$ is on $\S^{n-1}(c)$, then we say that $\S^{n-1}(c)$ or
$\B^n[c]$ \emph{supports} $C$ at $x$.
\end{defn}

The following corollary appears in \cite{KMP} without proof.
\begin{cor}\label{cor:intersection}
Let $A\subset\Re^n$ be a closed convex set. Then the following are equivalent.

\begin{tabular}{cl}
(i)& $A$ is spindle convex.\\
(ii)& $A$ is the intersection of unit balls containing it; that is, $A=\B[\B[A]]$.\\
(iii)& For every boundary point of $A$, there is a unit ball\\
& that supports $A$ at that point.
\end{tabular}
\end{cor}

\begin{thm}\label{thm:sep}
Let $C, D\subset\Re^n$ be spindle convex sets.
Suppose $C$ and $D$ have disjoint relative interiors. Then there is a closed unit ball $\B^n[c]$ 
such that $C\subseteq\B^n[c]$ and $D\subset\Re^n\backslash \B^n(c)$. 

Furthermore, if $C$ and $D$ have disjoint closures and one, say $C$, is
different from a  unit ball, then there is a closed unit ball $\B^n[c]$
such that $C \subset \B^n(c)$ and $D\subset\Re^n\backslash \B^n[c]$.
\end{thm}
 
\begin{proof}
Since $C$ and $D$ are spindle convex, they are convex, bounded sets
with disjoint relative interiors. So, their closures are convex, compact sets 
with disjoint relative interiors. 
Hence, they can be separated by a hyperplane $H$ that supports $C$ at 
a point, say $x$.
The closed unit ball $\B^n[c]$ of Lemma~\ref{lem:supportbyball} satisfies 
the conditions of the first statement. 

For the second statement, we assume that $C$ and $D$ have disjoint closures, so $\B^n[c]$ is 
disjoint from the closure of $D$ and remains so even after a sufficiently small translation. 
Furthermore, $C$ is a spindle convex set that is different from a unit ball, so 
$c \notin \conv(C \cap \S^{n-1}(c))$. Hence, there is a sufficiently small translation of 
$\B^n[c]$ that satisfies the second statement. 
\end{proof}

\begin{defn}
Let $C, D\subset\Re^n, c\in\Re^n, r>0$. We say that $\S^{n-1}(c,r)$
\emph{separates} $C$ from $D$ if $C\subseteq\B^n[c,r]$ and
$D\subseteq\Re^n\backslash \B^n(c,r)$, or $D\subseteq\B^n[c,r]$ and
$C\subseteq\Re^n\backslash \B^n(c,r)$. If $C\subseteq\B^n(c,r)$ and
$D\subseteq\Re^n\backslash \B^n[c,r]$, or $D\subseteq\B^n(c,r)$ and
$C\subseteq\Re^n\backslash \B^n[c,r]$, then we say that $C$ and $D$ are
\emph{strictly separated by} $\S^{n-1}(c,r)$.
\end{defn}

\section{A Kirchberger-type Theorem for Ball-polyhedra}\label{sec:Kirchberger}

The following theorem of Kirchberger is well known (e.g., \cite{Bar}).
If $A$ and $B$ are finite (resp. compact) sets in $\Re^n$ with the property
that for any set $T\subseteq A\cup B$
of cardinality at most $n+2$ the two sets $A\cap T$ and $B\cap T$
can be strictly separated by a hyperplane, then $A$ and $B$ can be strictly 
separated by a hyperplane.
We show that no similar statement holds for separation by unit spheres.

\begin{figure}[here]
\includegraphics[width=0.5\textwidth]{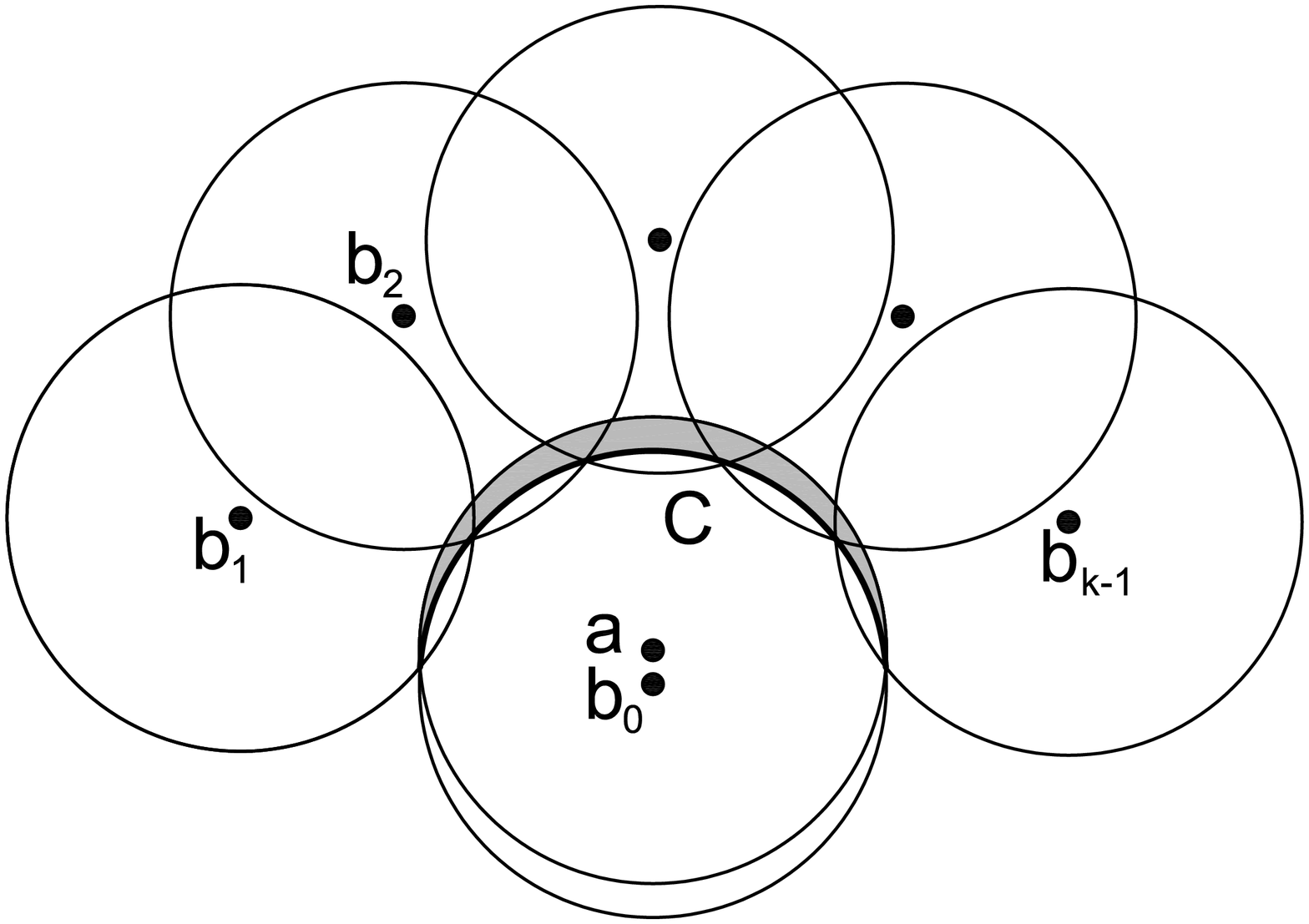}
\caption[]{}
\label{fig:kirchberger}
\end{figure}

We construct two sets $A$ and $B$ showing that there is no analogue of Kirchberger's theorem for 
separation by a unit sphere. Then we prove an analogue for separation by a sphere of radius at most one.
Let $A:=\{a\} \subset \Re^n$ be a singleton set and
$b_0 \in \Re^n$ be a point with $0 < ||a-b_0|| =: \delta < 1$.
Then $\B^n[a] \setminus \B^n(b_0)$ is a non-convex, closed set
bounded by two closed spherical caps: an inner one $C$ that belongs
to $\S^{n-1}(b_0)$ and an outer one that belongs to $\S^{n-1}(a)$
(Figure~\ref{fig:kirchberger}).
Now, we choose points $b_1, b_2, \ldots, b_{k-1}$ such that for every $i$ the set
$\B^n[b_i] \cap C$ is a spherical cap of radius $\varepsilon$
and we have also
\begin{equation}\label{eq:kirchberger}
C \subset \bigcup_{j=1}^{k-1} \B^n[b_j] \quad \hbox{and} \quad
C \not\subset \bigcup_{j=1, j\neq i}^{k-1} \B^n[b_j] \quad \hbox{for}\quad i=1,2,\ldots,k-1.
\end{equation}

\noindent
Let $B := \{ b_0, b_1, \ldots, b_{k-1} \}$. From (\ref{eq:kirchberger}) it easily follows that

\begin{equation}\label{eq:mainproperty}
\B^n(a) \subset \bigcup_{j=0}^{k-1} \B^n[b_j] \quad\hbox{and}\quad
\B^n(a) \not\subset \bigcup_{j=0, j\neq i}^{k-1} \B^n[b_j] \quad \hbox{for} \quad i=0,1,\ldots,k-1.
\end{equation}

From the first part of (\ref{eq:mainproperty}) it is clear that there is no $c \in \Re^n$
with the property that $a \in \B^n(c)$ and $B \subset \Re^n \setminus \B^n[c]$.
On the other hand, if $\varepsilon$ is sufficiently small, then $a \in \conv B$.
Hence, there is no $c \in \Re^n$ such that $B \subset \B^n(c)$ and $a \notin \B^n[c]$.
So, we have shown that $A$ and $B$ cannot be strictly separated by a unit sphere.

However, by the second part of (\ref{eq:mainproperty}),
for any $T \subset A \cup B$ of cardinality at most $k$, there is
a $c \in\Re^n$ such that $T \cap A \subset \B^n(c)$ and $T \cap B \subset \Re^n \setminus \B^n[c]$. 
This shows that \emph{there is no Kirchberger-type theorem for separation by unit spheres}.

In Theorem~\ref{thm:kirch_1conv} we provide a weaker analogue of Kirchberger's theorem.
For its proof we need the following version of Kirchberger's theorem,
which is a special case of Theorem 3.4 of Houle \cite{Ho91}, and a lemma.

\begin{thm}\label{thm:houle}
Let $A, B \subset \Re^n$ be finite sets.
Then $A$ and $B$ can be strictly separated by a sphere $\S^{n-1}(c,r)$ such that $A \subset \B^n(c,r)$
if, and only if, for every $T \subset A \cup B$ with $\card T  \leq n+2$,
$T \cap A$ and $T \cap B$ can be strictly separated by a sphere $\S^{n-1}(c_T,r_T)$ such that
$T \cap A \subset \B^n(c_T,r_T)$.
\end{thm}

\begin{lem}\label{lem:kirch_1conv}
Let $A, B \subset \Re^n$ be finite sets
and suppose that $\S^{n-1}(o)$ is the smallest sphere that separates $A$ from $B$ such that $A \subseteq \B^n[o]$.
Then there is a set $T \subseteq A \cup B$ with $\card T \leq n+1$
such that $\S^{n-1}(o)$ is the smallest sphere $\S^{n-1}(c,r)$ that separates
$T \cap A$ from $T \cap B$ and satisfies $T \cap A \subset \B^n[c,r]$.
\end{lem}

\begin{proof}
First observe that $A \neq \emptyset$.
Assume that $\S^{n-1}(o)$ separates $A$ from $B$ such that $A \subset \B^n[o]$.
Now, let us note also that $\S^{n-1}(o)$ is the smallest sphere separating $A$ and $B$ such that $A \subset \B^n[o]$
if, and only if, there is no closed spherical cap of radius less than $\pi /2$ that contains $A \cap \S^{n-1}(o)$ and
whose interior with respect to $\S^{n-1}(o)$ is disjoint from $B \cap \S^{n-1}(o)$.
Indeed, if there is a sphere $\S^{n-1}(x,r)$ that separates $A$ and $B$ and satisfies $r < 1$ and
$A \subset \B^n[x,r]$, then we may choose $\S^{n-1}(o) \cap \B^n[x,r]$ as such a spherical cap, a contradiction.
On the other hand, if $C$ is such a closed spherical cap then, by the finiteness of $A$ and $B$,
we can move $\S^{n-1}(o)$ to a sphere $\S^{n-1}(x,r)$ that separates $A$ and $B$
such that $\B^n[x,r] \cap \S^{n-1}(o) = C$ and $r < 1$, a contradiction.

We may assume that $A, B \subset \S^{n-1}(o)$. Let us take a point $q
\in \B^n[o]\setminus \{ o \}$. Observe that the closed half-space that
does not contain $o$ and whose boundary contains $q$ and is
perpendicular to $q$ intersects $\S^{n-1}(o)$ in a closed spherical cap
of radius less than $\pi /2$. Let us denote this spherical cap and its
interior with respect to $\S^{n-1}(o)$ by $C_q$ and $D_q$, respectively.
Observe that we have defined a one-to-one mapping between
$\B^n[o]\setminus \{ o\}$  and the family of closed spherical caps of
$\S^{n-1}(o)$ with radius less than $\pi /2$.

Let us consider a point $p \in \S^{n-1}(o)$.
Note that $p \in C_q$ for some $q \in \B^n[o]\setminus \{ o \}$ if, and only if,
the straight line passing through $p$ and $q$ intersects $\B^n[o]$ in a segment of
length at least $2 \| p-q \|$.

Set
\begin{equation}
F_p := \{ q \in \B^n[o]\setminus \{ o \} : p \in C_q \} \hbox{ and }
\end{equation}
$$
G_p := \{ q \in \B^n[o]\setminus \{ o \} : p \notin D_q \}.
$$

It is easy to see that 

\begin{equation}
F_p = \B^n[p/2,1/2] \setminus \{ o \} \hbox{ and }
G_p = \B^n[o] \setminus \Big( \B^n(p/2,1/2) \cup \{ o \} \Big).
\end{equation}

By the first paragraph of this proof,
$\S^{n-1}(o)$ is the smallest sphere separating $A$ and $B$ and satisfying $A \subset \B^n[o]$
if, and only if, $\big( \bigcap_{a \in A} F_a \big) \cap \big( \bigcap_{b \in B} G_b \big) = \emptyset$.

Let $f$ denote the inversion with respect to $\S^{n-1}(o)$; for the definition
of this transformation, we refer to \cite{Yag} Chapter III.
More specifically, let us define $f(x):= x / \| x \|^2$ for $x \in \Re^n \setminus \{ o \}$.
For any $p \in \S^{n-1}(o)$, let $H^+(p)$ (resp., $H^-(p)$), denote the closed half-space bounded 
by the hyperplane tangent to $\S^{n-1}(o)$ at $p$ that contains (resp., does not contain) $\S^{n-1}(o)$.
Using elementary properties of inversions, we see that $f(F_p) = H^-(p)$
and $f(G_p) = H^+(p) \setminus \B^n(o)$.
Hence, $\S^{n-1}(o)$ is the smallest sphere separating $A$ and $B$ and satisfying $A \subset \B^n[o]$ if,
and only if,

\begin{equation}\label{eq:helly}
I := \Big( \bigcap_{a \in A}  H^-(a) \Big) \cap \Big( \bigcap_{b \in B} \big( H^+(b)
\setminus \B^n(o)  \big) \Big)
\end{equation}

\noindent
is empty.
Observe that $\B^n(o) \cap H^-(a) = \emptyset$ for any $a \in A$.
Since $A \neq \emptyset$, we have
\begin{equation}
I = \Big( \bigcap_{a \in A}  H^-(a) \Big) \cap \Big( \bigcap_{b \in B} H^+(b) \Big).
\end{equation}

As $H^-(p)$ and $H^+(p)$ are convex for any $p \in \S^{n-1}(0)$, Helly's Theorem yields our statement.
\end{proof}

\begin{rem}
There are compact sets $A, B \subset \Re^n$ such that $\S^{n-1}(o)$ is the smallest sphere
that separates $A$ from $B$ and $A \subseteq \B^n[o]$ but, for any finite $T \subseteq A \cup B$,
there is a sphere $\S^{n-1}(x,r)$ that separates $T\cap A$ and $T\cap B$ such that $r < 1$ and
$T \cap A \subseteq \B^n[x,r]$.
\end{rem}

We show the following $3$-dimensional example.
Let us consider a circle $\S^1(x,r) \subset \S^2(o)$ with $r < 1$ and a set $A_0 \subset \S^1(x,r)$
that is the vertex set of a regular triangle.
Let $B$ be the image of $A_0$ under the reflection about $x$.
Clearly, $\S^1(x,r)$ is the only circle in its affine hull that separates $A_0$ and $B$.
Hence, every 2-sphere that separates $A_0$ and $B$ contains $\S^1(x,r)$.
Consider two points $a \in A_0$ and $y \in (o,a)$ and set $A = A_0 \cup \B^3(y, \| a-y \|)$.
Then the smallest sphere that separates $A$ and $B$ and contains $A$ in its convex hull is 
$\S^2(o)$.
Nevertheless, it is easy to show that, for any finite set $T \subset A$,
there is a sphere $\S^2(c_T,r_T)$ separating $T$ and $B$ such that $r_T < 1$ and $T \subset \B^3[c_T,r_T]$.

\begin{thm}\label{thm:kirch_1conv}
Let $A, B \subset \Re^n$ be finite sets. Then $A$ and $B$ can be
strictly separated by a sphere $\S^{n-1}(c,r)$ with $r\leq 1$ such that $A
\subset \B^n(c,r)$ if, and only if, the following holds. For every $T
\subseteq A \cup B$ with $ \card T \leq n+2$, $T \cap A$ and $T \cap B$
can be strictly separated by a sphere $\S^{n-1}(c_T,r_T)$ with $r_T\leq 1$
such that $T \cap A \subset \B^n(c_T,r_T)$.
\end{thm}

\begin{proof}
We prove the ``if'' part of the theorem, the opposite direction is trivial.
Theorem~\ref{thm:houle} guarantees the existence of the smallest 
sphere $\S^{n-1}(c',r')$
that separates $A$ and $B$ such that $A \subseteq \B^n[c',r']$.
According to Lemma~\ref{lem:kirch_1conv},
there is a set $T \subseteq A \cup B$ with $\card T \leq n+1$
such that $\S^{n-1}(c',r')$ is the smallest sphere that separates $T \cap A$
from $T \cap B$ and whose convex hull contains $T \cap A$.
By the assumption, we have $r' < r_T \leq 1$.
Note that Theorem~\ref{thm:houle} guarantees the existence of a sphere $\S^{n-1}(c^*,r^*)$
that strictly separates $A$ from $B$ and satisfies $A \subset \B^n(c^*,r^*)$.
Since $r' < 1$, there is a sphere $\S^{n-1}(c,r)$ with $r \leq 1$
such that $\B^n[c',r'] \cap \B^n(c^*,r^*) \subset \B^n(c,r) \subset \Re^n \setminus
\big( \B^n(c',r') \cup \B^n[c^*,r^*] \big)$.
This sphere clearly satisfies the conditions in Theorem~\ref{thm:kirch_1conv}.
\end{proof}

\begin{prob}
Prove or disprove that Theorem~\ref{thm:kirch_1conv} extends to compact sets.
\end{prob}

\section{The spindle convex hull:\\ The Theorems of Carath\'eodory and Steinitz}\label{sec:Caratheodory}
In this section we study the spindle convex hull of a set and give analogues of the
well-known theorems of Carath\'eodory and Steinitz to spindle convexity. 
The theorem of Carath\'eodory states that the convex hull of a set 
$X\subset\Re^n$ is the union of simplices with vertices in $X$.
Steinitz's theorem is that if a point is in the interior of the convex
hull of a set $X\subset\Re^n$, then it is also in the interior of the
convex hull of at most $2n$ points of $X$. This number $2n$ cannot be
reduced as shown by the cross-polytope and its center point. We state
the analogues of these two theorems in
Theorem~\ref{thm:CaratheodorySteinitz}. We note that, unlike in the case
of linear convexity, the analogue of the theorem of Kirchberger does not
imply the analogue of the theorem of Carath\'eodory.

Motivated by Lemma~\ref{lem:supportbyball} we make the following definition.

\begin{defn}
Let $X$ be a set in $\Re^n$. Then the \emph{spindle convex hull} of
$X$ is\\  
$\convv X := \bigcap \{ C\subseteq\Re^n : X \subseteq C$ 
and $C$ is spindle convex in $\Re^n\}$.
\end{defn}

The straightforward proof of the following elementary property of the spindle convex hull is omitted.

\begin{prop}
Let $P \subset H$, where $H$ is an affine subspace of $\Re^n$.
Assume that $A$ is contained in a closed unit ball.
Then the spindle convex hull of $P$ with respect to $H$ coincides
with the intersection of $H$ with the spindle convex hull of $P$
in $\Re^n$.
\end{prop}

\begin{defn}\label{defn:sphericallyconvex}
Let $\S^k(c,r)\subset\Re^n$ be a sphere such that $0\leq k \leq n-1$. A
set $F\subset \S^k(c,r)$ is \emph{spherically convex} if it is contained in an open
hemisphere of $\S^k(c,r)$ and for every $x,y\in F$ the shorter great-circular
arc of $\S^k(c,r)$ connecting $x$ with $y$ is in $F$. The \emph{spherical
convex hull} of a set $X\subset \S^k(c,r)$ is defined in the natural way and it
exists if, and only if, $X$ is in an open hemisphere of $\S^k(c,r)$. We denote
it by $\Sconv{X}{\S^k(c,r)}$.
\end{defn}

\begin{rem}\label{rem:sphericalcaratheodory}
Carath\'eodory's Theorem can be stated for the sphere in the following
way. If $X\subset\S^k(c,r)$ is a set in an open hemisphere of
$\S^k(c,r)$, then $\Sconv{X}{\S^k(c,r)}$ is the union of spherical
simplices with vertices in $X$. The proof of this spherical equivalent
of the theorem  uses the central projection of the open hemisphere to
$\Re^k$.
\end{rem}

\begin{rem}\label{rem:sphericallyconvexintersection}
It follows from Definition~\ref{defn:spindle} that if $C\subset\Re^n$ is a
spindle convex set such that $C\subset\B^n[q]$ and $\crr(C)<1$ then 
$C\cap\S^{n-1}(q)$ is spherically convex on $\S^{n-1}(q)$.
\end{rem}

\begin{figure}[here]
\includegraphics[width=0.33\textwidth]{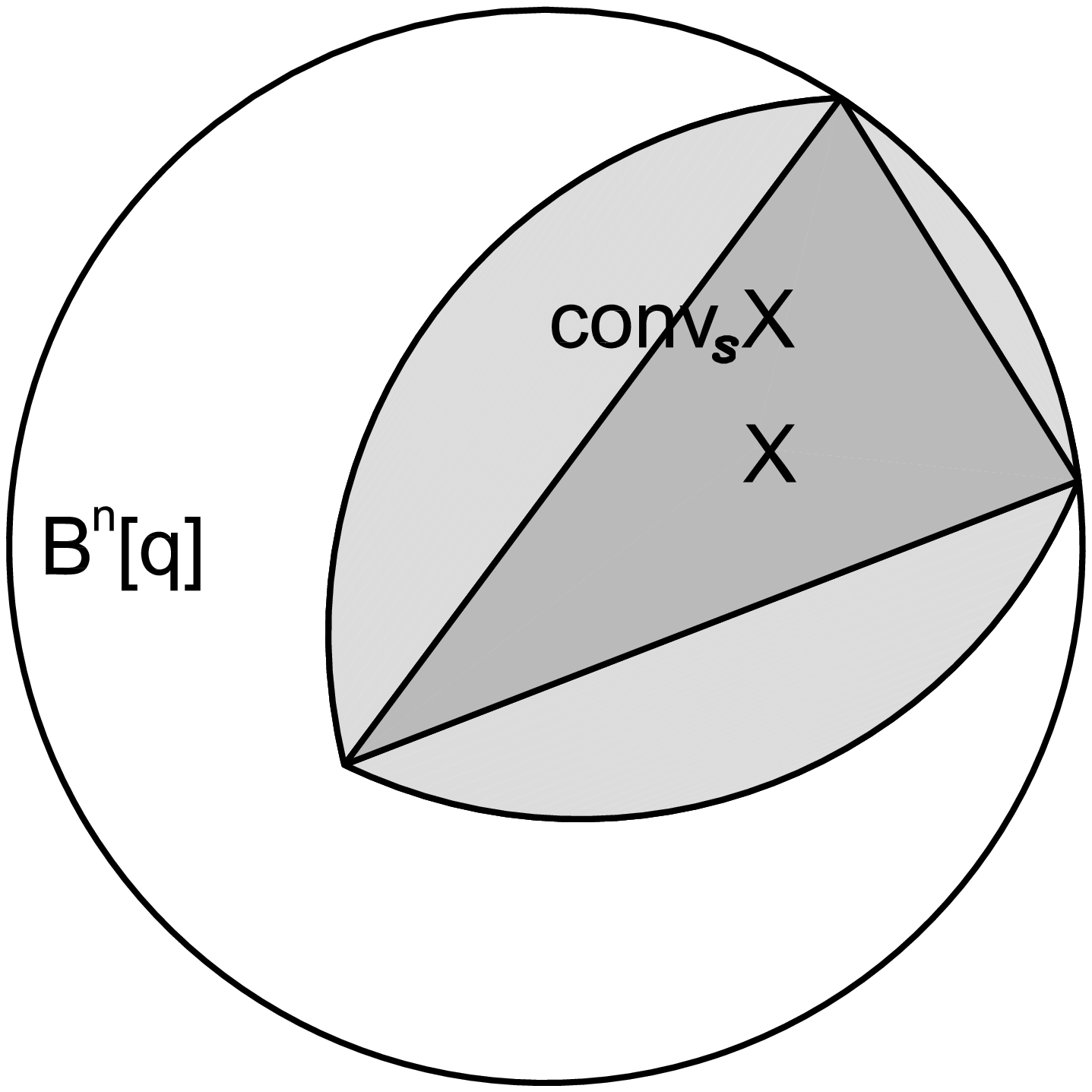}
\caption[]{}
\label{fig:caratheodory}
\end{figure}

The following lemma describes the surface of a spindle convex hull
(Figure~\ref{fig:caratheodory}).

\begin{lem}\label{lem:sconv}
Let $X\subset \Re^n$ be a closed set such that $\crr(X)<1$ and
let $\B^n[q]$ be a closed unit ball containing $X$.  Then

\begin{tabular}{cl}
(i)& 
$X\cap\S^{n-1}(q)$ is contained in an open hemisphere of $\S^{n-1}(q)$ and\\
(ii)&
$\convv(X)\cap\S^{n-1}(q)=\Sconv{X\cap\S^{n-1}(q)}{\S^{n-1}(q)}.$
\end{tabular}
\end{lem}

\begin{proof}
Since $\crr(X)<1$, we obtain that $X$ is contained in the intersection
of two distinct closed unit balls which proves (i). Note that by (i),
the right hand side $Z:=\Sconv{X\cap\S^{n-1}(q)}{\S^{n-1}(q)}$ of (ii)
exists.  We show that the set on the left hand side is contained in $Z$;
the other containment follows from 
Remark~\ref{rem:sphericallyconvexintersection}.

Suppose that $y\in\convv(X)\cap\S^{n-1}(q)$ is not contained in $Z$. We
show that there is a hyperplane $H$ through $q$ that strictly separates
$Z$ from $y$. Consider an open hemisphere of $\S^{n-1}(q)$ that contains
$Z$, call the spherical center of this hemisphere $p$. If $y$ is an
exterior point of the hemisphere, $H$ exists. If $y$ is on the boundary
of the hemisphere, then, by moving the hemisphere a little, we find
another open hemisphere that contains $Z$, but with respect to which $y$
is an exterior point.

Assume that $y$ is contained in the open hemisphere. Let $L$ be a
hyperplane tangent to $\S^{n-1}(q)$ at $p$. We project $Z$ and $y$
centrally from $q$ onto $L$ and, by the separation theorem of convex
sets in $L$, we obtain an $(n-2)$-dimensional affine subspace $T$ of $L$
that strictly separates the image of $Z$ from the image of $y$. Then
$H:=\aff(T\cup\{q\})$ is the desired hyperplane.

Hence, $y$ is contained in one open hemisphere of $\S^{n-1}(q)$ and $Z$
is in the other. Let $v$ be the unit normal vector of $H$ pointing
toward the hemisphere of $\S^{n-1}(q)$ that contains $Z$. Since $X$ is
closed, its distance from the closed hemisphere containing $y$ is
positive. Hence, we can move $q$ a little in the direction $v$ to obtain
the point $q'$ such that  $X\subset\B^n[q]\cap\B^n[q']$ and
$y\notin\B^n[q'].$  As $\B^n[q']$ separates $X$ from $y$, the latter is
not in $\convv X$, a contradiction.
\end{proof}

We prove the main result of this section.
\begin{thm}\label{thm:CaratheodorySteinitz}
Let $X\subset \Re^n$ be a closed set.

\begin{tabular}{cl}
(i)& If $y\in\bd\convv X$ then there is a set\\
&$\{x_1,x_2,\dots,x_n\}\subseteq X$ such that $y\in\convv \{x_1,x_2,\dots,x_n\}$.\\
(ii)& If $y\in\inter\convv X$ then there is a set\\
&$\{x_0,x_1,\dots,x_n\}\subseteq X$ such that $y\in\inter\convv \{x_0,x_1,\dots,x_n\}$.
\end{tabular}
\end{thm}

\begin{proof}
Assume that $\crr(X)>1$. Then $\B[X]=\emptyset$ hence, by Helly's
theorem, there is a set $\{x_0,x_1,\dots,x_n\}\subseteq X$
such that  $B[\{x_0,x_1,\dots,x_n\}]=\emptyset$. By
Corollary~\ref{cor:crr}, it follows that 
$\convv(\{x_0,x_1,\dots,x_n\})=\Re^n$. Thus, (i) and (ii) follow.

Now we prove (i) for $\crr(X)<1$. By
Lemma~\ref{lem:supportbyball}, Remark~\ref{rem:sphericalcaratheodory}
and Lemma~\ref{lem:sconv} we obtain that
$y\in\Sconv{\{x_1,x_2,\dots,x_n\}}{\S^{n-1}(q)}$ for some
$\{x_1,x_2,\dots,x_n\}\subset X$ and some $q\in\Re^n$ such that
$X\subseteq\B^n[q]$.  Hence, $y\in\convv \{x_1,x_2,\dots,x_n\}$.

We prove (i) for $\crr(X)=1$ by a limit argument as follows. Without
loss of generality, we may assume that $X\subseteq\B^n[o]$. Let
$X^k:=(1-\frac{1}{k})X$ for any $k\in\Ze^+$. Let $y^k$ be the point of
$\bd \convv (X^k)$ closest to $y$. Thus,
$\lim\limits_{k\rightarrow\infty}y^k=y$. Clearly, $\crr(X^k)<1$, hence there is
a set $\{x^k_1,x^k_2,\dots,x^k_n\}\subseteq X^k$ such that
$y^k\in\convv\{x^k_1,x^k_2,\dots,x^k_n\}$. By compactness,
there is a sequence $0<i_1<i_2<\dots$ of indices such that all the $n$
sequences  
$\{x^{i_j}_1:j\in\Ze^+\},\{x^{i_j}_2:j\in\Ze^+\},\dots,\{x^{i_j}_n:j\in\Ze^+\}$
converge. Let their respective limits be $x_1,x_2,\dots,x_n$. Since $X$
is closed, these $n$ points are contained in $X$. Clearly,  $y\in\convv
\{x_1,x_2,\dots,x_n\}$.

To prove (ii) for $\crr(X)\leq 1$, suppose that
$y\in\inter\convv X$. Then let $x_0\in X\cap\bd\convv X$ be arbitrary
and let $y_1$ be the intersection of $\bd\convv X$ with the ray starting
from $x_0$ and passing through $y$. Now, by (i), $y_1\in\convv
\{x_1,x_2,\dots,x_n\}$ for some $\{x_1,x_2,\dots,x_n\}\subseteq X$.
Then clearly $y\in\inter\convv
\{x_0,x_1,\dots,x_n\}$.
\end{proof}

The same proof with a simple modification provides the analogue of
the ``Colorful Carath\'eodory Theorem'' 
(\cite{Matousek} p. 199).

\begin{thm}
Consider $n+1$ finite point sets $X_1,\dots,X_{n+1}$ in
$\Re^n$ such that the spindle convex hull of each contains the origin.
Then there is an $(n+1)$-point set $T\subset X_1\cup\dots\cup X_{n+1}$ 
with $\card(T\cap X_i)=1$ for each $i\in\{1,2,\dots,n+1\}$ 
such that $o\in\convv T$.
\end{thm}

\section{The Euler-Poincar\'e Formula for Standard Ball-polyhedra}\label{sec:EPF}

The main result of this section is the Euler-Poincar\'e formula for a
certain family of ball-polyhedra. However, before developing that, we
present Example~\ref{ex:barrel} to show that describing the face
lattice of arbitrary ball-polyhedra is a difficult task. The example
is as follows.

We construct a $4$-dimensional ball-polyhedron $P$ which has a subset $F$
on its boundary that, according to any meaningful definition of a face
for ball-polyhedra, is a $2$-dimensional face. However, $F$ is
homeomorphic to a band, hence it is \emph{not} homeomorphic to a
disk.  This example demonstrates that even if one finds a satisfactory
definition for the face lattice of a ball-polyhedron that models the face
lattice of a convex polytope, it will not lead to a CW-decomposition of
the boundary of ball-polyhedra. 

\begin{ex}\label{ex:barrel} Take two unit spheres in $\Re^4$, $\S^3(p)$ and
$\S^3(-p)$ that intersect in a $2$-sphere
$\S^2(o,r):=\S^3(p)\cap\S^3(-p)$ of $\Re^4$. Now, take a closed unit
ball $\B^4[q]\subset\Re^4$ that intersects $\S^2(o,r)$ in a spherical
cap $\S^2(o,r)\cap\B^4[q]$ of $\S^2(o,r)$ which is greater than a
hemisphere of $\S^2(o,r)$, but is \emph{not} $\S^2(o,r)$. Such a unit ball
exists, since $r<1$. Let the ball-polyhedron be
$P:=\B^4[p]\cap\B^4[-p]\cap\B^4[q]\cap\B^4[-q]$. Now,
$F:=\S^2(o,r)\cap\B^4[q]\cap\B^4[-q]$ is homeomorphic to a
two-dimensional band. Also, $F$ is a subset of the boundary of $P$ that
deserves the name of ``2-face''.
\end{ex}

\begin{defn}
Let $\S^l(p,r)$ be a sphere of $\Re^n$. The intersection of $\S^l(p,r)$
with an affine subspace of $\Re^n$ that passes through $p$ is called a
\emph{great-sphere} of $\S^l(p,r)$. Note that $\S^l(p,r)$ is a great-sphere of
itself. Moreover, any great-sphere is itself a sphere.
\end{defn}

{\defn Let $P\subset\Re^n$ be a ball-polyhedron with a family of
generating balls $\B^n[x_1],\dots,\B^n[x_k]$. This family of generating
balls is called \emph{reduced} if removing any of the balls yields that
the intersection of the remaining balls becomes a set larger than $P$. 
Note that, for any ball-polyhedron, distinct from a singleton, 
there is a unique reduced family of generating balls. 
A \emph{supporting sphere} $\S^l(p,r)$ of $P$ is a
sphere of dimension $l$, where $0\leq l\leq (n-1)$, which can be
obtained as an intersection of some of the generating spheres of $P$
from the reduced family of generating spheres of $P$ such that
$P\cap\S^l(p,r)\neq\emptyset$.}

Note that the intersection of finitely many spheres in $\Re^n$ is either
empty, or a sphere, or a point.

In the same way that the faces of a convex polytope can be described in terms
of supporting hyperplanes, we describe the faces of a certain class of
ball-polyhedra in terms of supporting spheres. 

{\defn  Let $P$ be an $n$-dimensional ball-polyhedron. We say that $P$ is
\emph{standard} if for any supporting sphere
$\S^l(p,r)$ of $P$ the intersection $F:=P \cap \S^l(p,r)$ is
homeomorphic to a closed Euclidean ball of some dimension. We call $F$ a
\emph{face} of $P$, the \emph{dimension} of $F$ is the dimension of the ball
that $F$ is homeomorphic to. If the dimension is $0, 1$ or $n-1$, then
we call the face a \emph{vertex}, an \emph{edge} or a \emph{facet},
respectively.}

Note that the dimension of $F$ is independent of the choice of the
supporting sphere containing $F$.

In Section~\ref{sec:Steinitz}, we present reasons why standard
ball-polyhedra are natural, relevant objects of study in $\Re^3$. 
Example~\ref{ex:barrel} demonstrates the reason behind studying these
objects in higher dimensions.

%\begin{claim}
%Let $P$ be a standard ball-polyhedron with the reduced family of
%generating balls, $B[p_1],\dots,\B[p_k]$, in $\Re^n$. Then
%$\aff\{p_1,\dots,p_k\}=\Re^n$.
%\end{claim}
%\begin{proof}
%Let $H:=\aff\{p_1,\dots,p_k\}$ and assume that $l:=\dim H<n$. We will
%show that $P$ is \emph{not} standard. Let $q\in H$ be the center of the
%smallest ball (the circumball of the points) in $H$ that contains
%$\{p_1,\dots,p_k\}$. Since  $P \cap H$ has
%non-empty relative interior in $H$, the radius of the circumball is less
%than $1$. Now, let $q'\in\Re^n$ be a point such that $\S^{n-1}(q')$
%contains the boundary of the circumball. So, $B[q]\supset
%\{p_1,\dots,p_k\}$. Clearly, $\S(q)$ contains at least $l+1$ points from
%$\{p_1,\dots,p_k\}$, and those that it contains affinely span $H$. We
%select $l+1$ affinely independent points of these points, and by the
%appropriate choice of indeces we can assume that the selected ones are
%$\{p_1,\dots,p_{l+1}\}$. Now,
%$q'\in\cap\{\S^{n-1}(p_1),\dots,\S^{n-1}(p_{l+1})\}$ and $q'\in P$.
%Moreover, if we arbitrarily rotate $q'$ about $H$ (i.e. we take any
%isometry of $\Re^n$ which acts as the identity on $H$) to get $q''$ we have
%that $q''\in P$. The orbit of $q''$ under all rotations about $H$ is a
%sphere, $\S^{n-l-1}(c,r)$. By the choice of $p_1,\dots, p_{l+1}$ we have
%that $\S^{n-l-1}(c,r)=\cap\{\S^{n-1}(p_1),\dots,\S^{n-1}(p_{l+1})\}$.
%Since $q\in \bd P$ and $P$ is invariant under rotations about $H$, we
%have $\bd P \supset \S^{n-l-1}(c,r)$. Hence, $P$ is \emph{not} standard.
%\end{proof}

For the proof of the next theorem we need the following definition.

\begin{defn}
Let $K$ be a convex body in $\Re^n$ and $b\in\bd K$. Then the
\emph{Gauss image} of $b$ with respect to $K$ is the set of outward unit
normal vectors of hyperplanes that support $K$ at $b$. Clearly, it is a
spherically convex subset of $\S^{n-1}(o)$ and its dimension is defined
in the natural way.
\end{defn}

\begin{thm}\label{thm:CW}
Let $P$ be a standard ball-polyhedron. Then the faces of $P$ form
the closed cells of a finite CW-decomposition of the boundary of $P$.
\end{thm}
\begin{proof}
Let $\{\S^{n-1}(p_1),\dots,\S^{n-1}(p_k)\}$ be the reduced family of
generating spheres of  $P$. The \emph{relative interior} (resp., the
\emph{relative boundary}) of an $m$-dimensional face $F$ of $P$ is
defined as the set of those points of $F$ that are mapped to $\B^m(o)$
(resp., $\S^{m-1}(o)$) under any homeomorphism between $F$ and
$\B^m[o]$.

For every $b\in\bd P$ define the following sphere
\begin{equation}
S(b):=\bigcap\{\S^{n-1}(p_i) : p_i\in\S^{n-1}(b), i\in\{1,\dots,k\}\}.
\end{equation}
Clearly, $S(b)$ is a support sphere of $P$. Moreover, if $S(b)$ is an
$m$-dimensional sphere, then the face $F:=S(b)\cap P$ is also
$m$-dimensional as $b$ has an $m$-dimensional neighbourhood in $S(b)$
that is contained in $F$. This also shows that $b$ belongs to the
relative interior of $F$. Hence, the union of the relative interiors of
the faces covers $\bd P$.

We claim that every face $F$ of $P$ can be obtained in this way, i.e., for
any relative interior point $b$ of $F$ we have $F=S(b)\cap P$.
Clearly, $F\supseteq S(b)\cap P$, as the support sphere of
$P$ that intersects $P$ in $F$ contains $S(b)$. It is sufficient to show  
that $F$ is at most $m$-dimensional. This is so, because the Gauss image
of $b$ with respect to $P$ is at least $(n-m-1)$-dimensional, since the
Gauss image of $b$ with respect to 
$\bigcap\{\B^{n}[p_i] : p_i\in\S^{n-1}(b), i\in\{1,\dots,k\}\}\supseteq P$
is $(n-m-1)$-dimensional.

The above argument also shows that no point $b\in \bd P$ belongs to
the relative interior of more than one face. Moreover, if $b\in \bd P$
is on the relative boundary of the face $F$ then $S(b)$ is clearly of
smaller dimension than $F$. Hence, $b$ belongs to the relative interior
of a face of smaller dimension. This concludes the proof of the theorem.
\end{proof}

\begin{cor}
The reduced family of generating balls of any standard ball-polyhedron
$P$ in $\Re^n$ consists of at least $n+1$ unit balls.
\end{cor}

\begin{proof}
Since the faces form a CW-decomposition of the boundary of $P$, it has a
vertex $v$. The Gauss image of $v$ is $(n-1)$-dimensional. So, $v$
belongs to at least $n$ generating spheres from a reduced family. We
denote the centers of those spheres by $x_1, x_2,\dots,x_n$. Let
$H:=\aff\{x_1, x_2,\dots,x_n\}$. Then $\B[\{x_1, x_2,\dots,x_n\}]$ is
symmetric about $H$. Let $\sigma_H$ be the reflection of $\Re^n$ about
$H$. Then $S:=\S(x_1)\cap\S(x_2)\cap\dots\cap\S(x_n)$ contains the
points $v$ and $\sigma_H(v)$, hence $S$ is a sphere, not a point. Since
$P$ is a standard ball-polyhedron, there is a unit-ball
$\B[x_{n+1}]$ in the reduced family of generating balls of $P$ that does
not contain $S$.
\end{proof}

\begin{cor}
Let $\Lambda$ be the set containing all faces of a standard ball-polyhedron
$P\subset\Re^n$ and the empty set and $P$ itself. Then $\Lambda$ is a finite
bounded lattice with respect to ordering by inclusion. The atoms of
$\Lambda$ are the vertices of $P$ and $\Lambda$ is atomic, i.e., for every
element $a\in\Lambda$ with $a\neq\emptyset$ there is a vertex $x$ of
$P$ such that $x\in a$.
\end{cor}
\begin{proof}
First, we show that the intersection of two faces $F_1$ and $F_2$ is
another face (or the empty set). The intersection of the two supporting
spheres that intersect $P$ in $F_1$ and $F_2$ is another supporting
sphere of $P$, say $\S^l(p,r)$. Then $\S^l(p,r)\cap P=F_1\cap F_2$ is a
face of $P$. From this the existence of a unique maximum common lower
bound (i.e. an \emph{infimum}) for $F_1$ and $F_2$ follows. 

Moreover, by the finiteness of $\Lambda$, the existence of a
unique infimum for any two elements of $\Lambda$ implies the existence
of a unique minimum common upper bound (i.e., a \emph{supremum}) for any two
elements of $\Lambda$, say $C$ and $D$, as follows. The supremum of $C$
and $D$ is the infimum of all the (finitely many) elements of $\Lambda$
that are above $C$ and $D$. 

Vertices of $P$ are clearly atoms of $\Lambda$. Using
Theorem~\ref{thm:CW} and induction on the dimension of the face it is
easy to show that every face is the supremum of its vertices. 
\end{proof}

\begin{cor}\label{cor:existenceoffaces}
A standard ball-polyhedron $P$ in $\Re^n$ has $k$-dimensional faces for every
$0\leq k\leq n-1$.
\end{cor}
\begin{proof}
We use an inductive argument on $k$, where we go from $k=n-1$ down to
$k=0$. Clearly, $P$ has facets. A $k$-face $F$ of $P$ is homeomorphic
to $\B^k[o]$, hence its relative boundary is homeomorphic to $\S^{k-1}$,
if $k>0$. Since the $(k-1)$-skeleton of $P$ covers the relative boundary
of $F$, $P$ has $(k-1)$-faces.
\end{proof}

{\cor (Euler-Poincar\'e Formula) For any standard $n$-dimensional
ball-polyhedron $P$ we have:
\[
1 + (-1)^{n+1} = \sum_{i=0}^{n-1} (-1)^i f_i(P), 
\]
where $f_i(P)$ denotes the number of $i$-dimensional faces
of $P$.} 
\begin{proof}
It follows from the above theorem and the fact
that a ball-polyhedron in $\Re^n$ is a convex body, hence
its boundary is homeomorphic to $\S^{n-1}(o)$.
\end{proof}

\begin{cor}
Let $n\geq 3$. Any standard ball-polyhedron $P$ is the spindle convex hull
of its $(n-2)$-dimensional faces. Furthermore, no standard ball-polyhedron
is the spindle convex hull of its $(n-3)$-dimensional faces.
\end{cor}
\begin{proof}
For the first statement, it is sufficient to show that the spindle convex
hull of the $(n-2)$-faces contains the facets. Let $p$ be a point on the
facet, $F=P\cap\S^{n-1}(q)$. Take any great circle $C$ of
$\S^{n-1}(q)$. Since $F$ is spherically convex on $\S^{n-1}(q)$, $C\cap
F$ is a unit circular arc of length less than $\pi$. Let
$r,s\in\S^{n-1}(q)$ be the two endpoints of $C\cap F$. Then $r$ and $s$
belong to the relative boundary of $F$. Hence, by Theorem~\ref{thm:CW},
$r$ and $s$ belong to an $(n-2)$-face. Clearly, $p\in\convv \{r,s\}$.

The proof of the second statement follows. By
Corollary~\ref{cor:existenceoffaces} we can choose a relative interior
point $p$ of an $(n-2)$-dimensional face $F$ of $P$. Let $q_1$ and
$q_2$ be the centers of the generating balls of $P$ from a reduced
family such that $F:=\S^{n-1}(q_1)\cap\S^{n-1}(q_2)\cap P$. Clearly,
$p\notin\convv((\B[q_1]\cap\B[q_2])\backslash\{p\})\supseteq
\convv(P\backslash\{p\})$.
\end{proof}

\section{A Counterexample to a Conjecture of Maehara\\ in Dimensions at Least Four}\label{sec:Maehara}

Helly's Theorem, as stated for convex sets, adds nothing to the current theory.
However, the following result of Maehara \cite{Mae} is very suggestive.

\begin{thm}\label{thm:maehara}
Let $\F$ be a family of at least $n+3$ distinct $(n-1)$-spheres in $\Re^n$.
If any $n+1$ of the spheres in $\F$ have a point in common,
then all of the spheres in $\F$ have a point in common.
\end{thm}

Maehara points out that neither $n+3$ nor $n+1$ can be reduced.
First, we prove a variant of Theorem~\ref{thm:maehara}.

\begin{thm}
Let $\mathfrak{F}$ be a family of $(n-1)$-spheres in $\Re^n$, and $k$ be an 
integer such that $0 \leq k \leq n-1$. Suppose that $\mathfrak{F}$ has at 
least $n-k$ members and that any $n-k$ of them intersect in a sphere of 
dimension at least $k+1$. Then they all intersect in a sphere of 
dimension at least $k+1$. Furthermore, $k+1$ cannot be reduced to $k$.
\end{thm}

\begin{proof}
Amongst all the intersections of any $n-k$ spheres from the family, let $S$ be
such an intersection of minimal dimension. By assumption, $S$ is a sphere of
dimension at least $k+1$. Now, one of the $n-k$ spheres is redundant in the sense
that $S$ is contained entirely in this sphere. After discarding this redundant sphere, $S$ is now the
intersection of only $(n-k)-1$ members of the family, but any $n-k$ members
intersect in a sphere of dimension at least $k+1$. So, the remaining members of the
family intersect $S$. Since the dimension of $S$ is minimal, $S$ is contained in these members.
In particular, $\bigcap \mathfrak{F} = S$.

Fixing $n$ and $k$, $0 \leq k \leq n-1$, we show that $k+1$ cannot be
reduced to $k$ by considering a regular $n$-simplex in $\Re^n$, with
circumradius one, and a family of $n+1$ unit spheres centered at the
vertices of this simplex. The intersection of any $n-k$ of them is a
sphere of dimension at least $k$, but the intersection of all of them is a
single point which, as we recall, is not a sphere in the current setting.
\end{proof}

Maehara \cite{Mae} conjectured the following stronger version of
Theorem~\ref{thm:maehara}.

\begin{conj}\label{conj:Maehara}
Let $\F$ be a family of at least $n+2$ distinct $(n-1)$-dimensional unit spheres in $\Re^n$,
where $n\geq 3$. Suppose that any $n+1$ spheres in $\F$ have a point 
in common. Then all the spheres in $\F$ have a point in common.
\end{conj}

After Proposition~3 in \cite{Mae}, Maehara points out the importance
of the condition $n \geq 3$ by showing the following statement, also known as 
\emph{\c Ti\c teica's theorem} (sometimes called \emph{Johnson's theorem}). 
This theorem was found by the Romanian mathematician
G. \c Ti\c teica in 1908 (for historical details, see also \cite{BB}, and \cite{Johns}, p. 75).

\begin{figure}[here]
\includegraphics[width=0.37\textwidth]{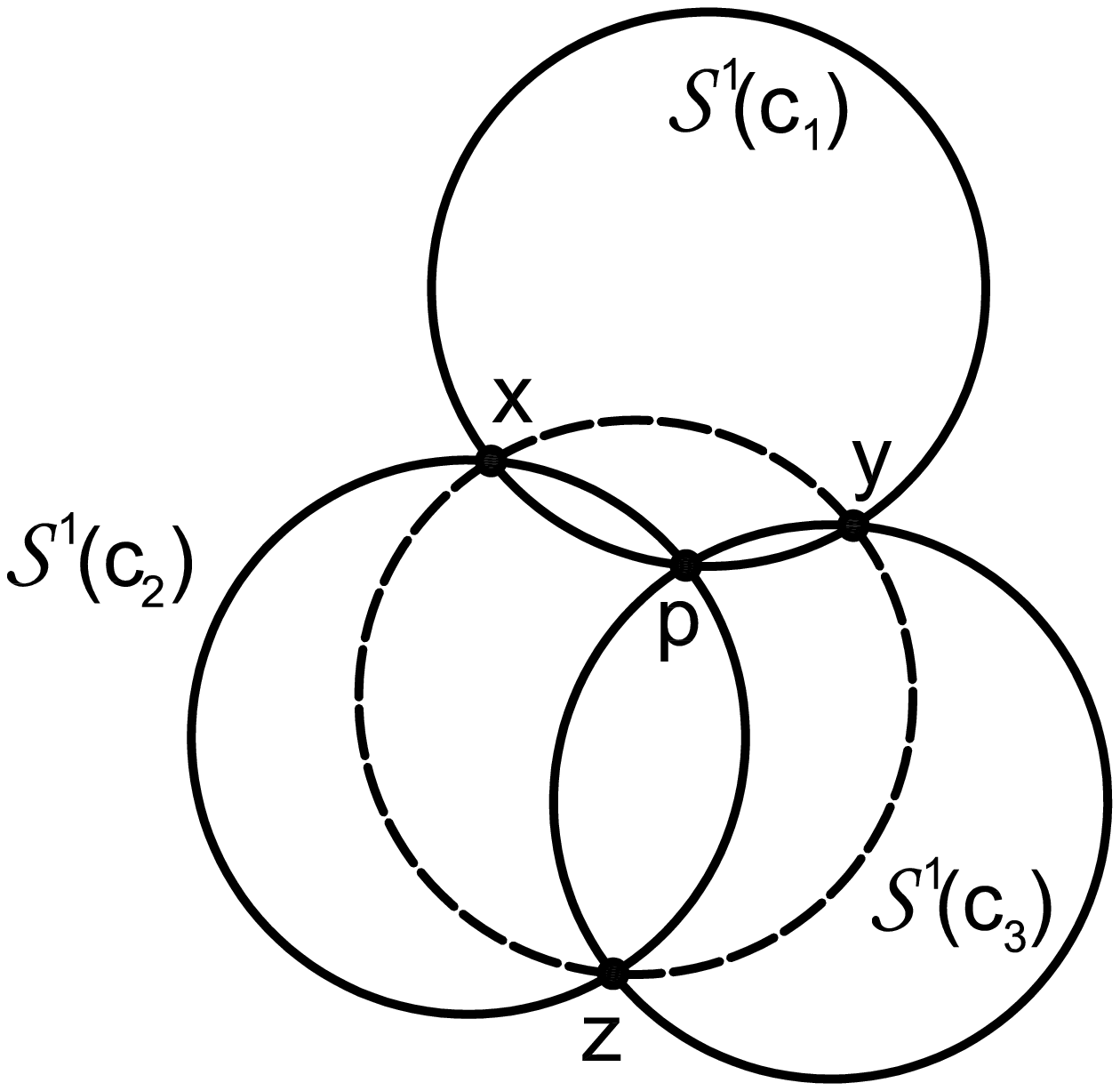}
\caption[]{}
\label{fig:4circle}
\end{figure}

\begin{prop}\label{prop:4circle}
Let $\S^1(c_1)$, $\S^1(c_2)$ and $\S^1(c_3)$ be unit circles in $\Re^2$
that intersect in a point $p$ (see Figure~\ref{fig:4circle}).
Let $\{x, p\}:= \S^1(c_1) \cap \S^1(c_2)$, $\{y, p\}:= \S^1(c_1) \cap \S^1(c_3)$ and
$\{z, p\}:= \S^1(c_2) \cap \S^1(c_3)$.
Then $x$, $y$ and $z$ lie on a unit circle.
\end{prop}

In the remaining part of this section, we show that
Conjecture~\ref{conj:Maehara} is false for $n \geq 4$.
To construct a suitable family $\F$ of unit spheres, we need the following lemma.

\begin{lem}\label{lem:maehara counterexample}
The following are equivalent.

\begin{tabular}{cl}
(i) & There is an $n$-simplex $P \subset \Re^n$ with circumsphere $\S^{n-1}(o,R)$ and\\
& a sphere $\S^{n-1}(x_1,r)$ tangent to all facet-hyperplanes of $P$ such that\\
& either $R^2 - 2rR =d^2$ or $R^2 + 2rR = d^2$ holds, where
$d:=\| x_1 - o\|$.\\
(ii) & There is a family of $n+2$ distinct $(n-1)$-dimensional unit spheres in\\
&  $\Re^n$ such that any $n+1$ of them have a common point but not all of\\
& them have a common point.
\end{tabular}
\end{lem}

\begin{proof}
First, we show that $(ii)$ follows from $(i)$.
Observe that, from $R^2 - 2rR = d^2$, we have
$R > d$, which implies that $x_1 \in \B^n(o,R)$.
Similarly, if $R^2 + 2rR = d^2$, then $x_1 \notin \B^n[o,R]$.
Thus, $x_1 \notin \S^{n-1}(o,R)$.
Since $\S^{n-1}(x_1,r)$ is tangent to every facet-hyperplane of $P$, $x_1$ is not
contained in any of these hyperplanes.

Consider the inversion $f$ with respect to $\S^{n-1}(x_1,r)$.
Let $a_i$ be a vertex of $P$ and $H_i$ denote the facet-hyperplane of $P$
that does not contain $a_i$, for $i=2, 3, \ldots, n+2$.
We set $\S^{n-1}(c_i,r_i):=f(H_i)$, $x_i:=f(a_i)$, for $i=2,3,\ldots,n+2$.
Finally, $\S^{n-1}(c_1,r_1):= f(\S^{n-1}(o,R))$.

Let $2 \leq i \leq n+2$.
Since $H_i$ is tangent to $\S^{n-1}(x_1,r)$, $S^{n-1}(c_i,r_i)$
is a sphere tangent to $\S^{n-1}(x_1,r)$ and contains $x_1$.
Hence, the radius of $\S^{n-1}(c_i,r_i)$ is $r_i = \frac{r}{2}$.
We show that also $r_1 = \frac{r}{2}$.
If $x_1 \in \B^n(o,R)$, then, using the definition of inversion and the equations
in  $(i)$, we have
\begin{equation}
2r_1 = \diam \S^{n-1}(c_1,r_1) = \frac{r^2}{R+d} + \frac{r^2}{R-d} = \frac{2r^2R}{R^2-d^2} = r.
\end{equation}
If $x_1 \notin \B^n[o,R]$, then

\begin{equation}
2r_1 = \diam \S^{n-1}(c_1,r_1) = \frac{r^2}{d-R} - \frac{r^2}{d+R} = \frac{2r^2R}{d^2-R^2} = r.
\end{equation}

Let $\F := \{ \S^{n-1}(c_i,\frac{r}{2}) : i= 1, \ldots, n+2 \}$.
Observe that $x_1 \in \S^{n-1}(c_i,\frac{r}{2})$, for every $i \neq 1$, and that
$x_i \in \S^{n-1}(c_1,\frac{r}{2}) \cap \S^{n-1}(c_j,\frac{r}{2})$ for every $j \neq i$.
So, $\F$ is a family of $n+2$ spheres of radius $\frac{r}{2}$ such that any $n+1$ of them
have a common point.

We assume $y \in \bigcap \F$.
Then $y \neq x_1$, since $x_1 \notin \S^{n-1}(o,R)$ and thus
$x_1 \notin \S^{n-1}(c_1,\frac{r}{2}) = f(\S^{n-1}(o,R))$.
Hence, $z := f(y) = f^{-1}(y)$ exists.
The point $z$, according to our assumption, is contained in every facet-hyperplane of $P$ and
also in its circumsphere, a contradiction.
So, $\F' := \{ \S^{n-1}(\frac{2}{r}\cdot c_i) : i= 1, \ldots, n+2 \}$ is a family of
unit sphere that satisfies $(ii)$.

A similar argument shows that $(ii)$ implies $(i)$.
\end{proof}

\begin{thm}\label{thm:counterMaehara}
For any $n \geq 4$, there exists a family of $n+2$ distinct $(n-1)$-dimensional unit spheres
in $\Re^n$ such that any $n+1$, but not all, of them have a common point.
\end{thm}

\begin{proof}
We use Lemma~\ref{lem:maehara counterexample} and construct a simplex $P$ and a sphere
$\S^{n-1}(x_1,r)$ such that they satisfy Lemma~\ref{lem:maehara counterexample} $(i)$.
We set $m:=n-1$.

Consider a line $L$ containing the origin $o$ and a hyperplane $H$
which is orthogonal to $L$ and is at a given distance $t \in (0,1)$ from $o$.
Let $u$ denote the intersection point of $L$ and $H$.
We observe that $t = \| u \|$ and let $b:= \frac{1}{t}u$.
Then $b \in \S^{n-1}(o,1)$.
Let $F$ be a regular $m$-simplex in $H$ whose circumsphere
is $S^{n-1}(o,1) \cap H$.
Thus, $u$ is the center of $F$ and the circumsphere of
$P:=\conv(F \cup \{b\})$ is $\S^{n-1}(o,1)$.
Clearly, there is a unique sphere $\S^{n-1}(c,r)$
tangent to every facet-hyperplane of $P$ such that
$c \in L$ and $c \notin P$.
We set $d:=\|c\|$.

Our aim is to show that, with a suitable choice of $t$,
$P$ and $\S^{n-1}(c,r)$ satisfy Lemma~\ref{lem:maehara counterexample} $(i)$.
So, in the remaining part of the proof, we calculate $g_m(t):=d(t)^2+2r(t)-1$ as a function of $t$,
and show that this function has a root on the interval $(0,1)$, for $m \geq 3$.
We note that if $g_m(t) = 0$, for some $t$, then $P$ and $\S^{n-1}(c,r)$ satisfy
the first equality in Lemma~\ref{lem:maehara counterexample} $(i)$ for $R = 1$
and $x_1=c$.

\begin{figure}[here]
\includegraphics[width=0.42\textwidth]{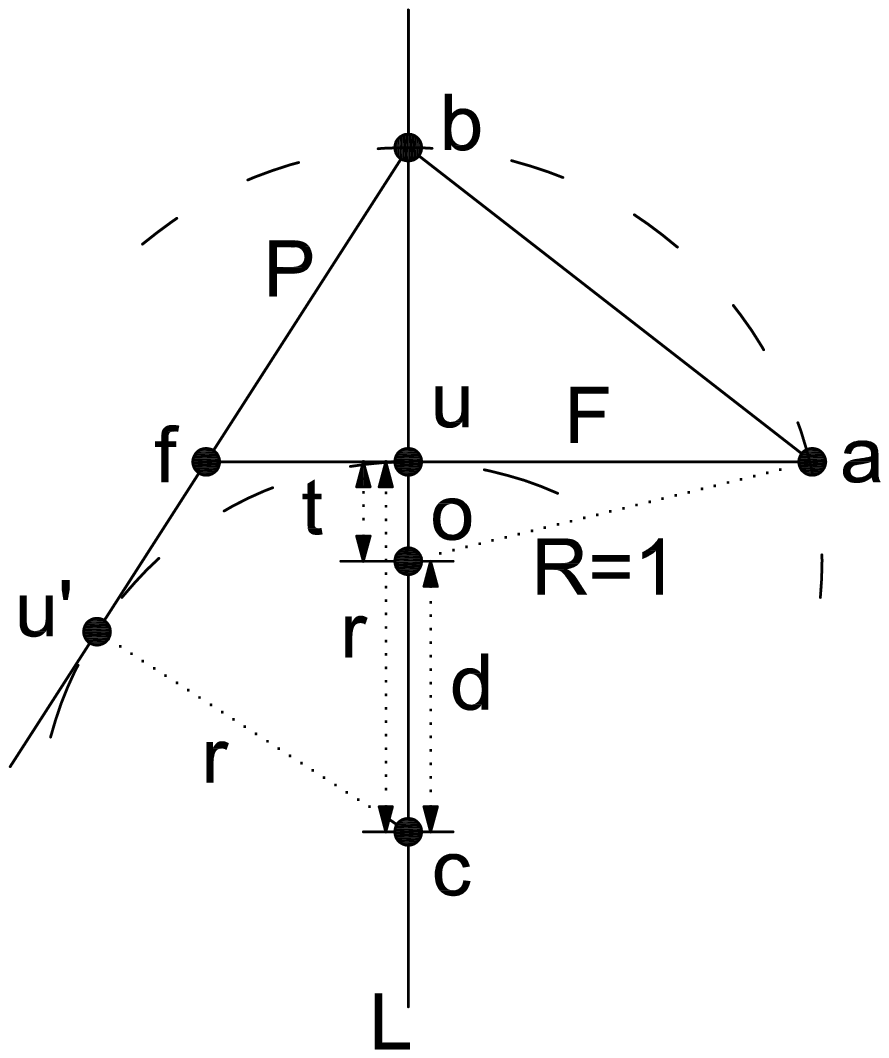}
\caption[]{}
\label{fig:maehara}
\end{figure}

Consider a vertex $a$ of $F$ and the center $f$ of the facet of $F$
that does not contain $a$ (see Figure~\ref{fig:maehara}).
Then $\|b-u\| = 1-t$ and $\|a-u\| = \sqrt{1-t^2}$.
Since in an $m$-dimensional regular simplex the distance of the center from any 
vertex is $m$ times
as large as the distance of the center and any facet-hyperplane, we have
$\|u-f\| = \frac{\sqrt{1-t^2}}{m}$.
We observe that $\S^{n-1}(c,r)$ is tangent to the facet-hyperplane $H_a$ of $P$ that
does not contain $a$.
Let $u'$ denote the intersection point of $\S^{n-1}(c,r)$ and $H_a$.
Clearly, $u'$, $f$ and $b$ are collinear and $\|u-c\| = \|u'-c\| = r$.
Furthermore, the two triangles $\conv \{ c,u',b \}$ and $\conv \{ f,u,b \}$ are co-planar and similar.
Hence,

\begin{equation}
\frac{\|b-f\|}{\|b-c\|} = \frac{\|u-f\|}{\|u'-c\|}.
\end{equation}

We have that $\|b-f\| = \sqrt{(1-t)^2+\frac{1-t^2}{m^2}}$, $\|b-c\| = 1+r-t$,
$\|u'-c\| = r$ and $\|u-f\| = \frac{\sqrt{1-t^2}}{m}$.
So, we have an equation for $r$ which yields

\begin{equation}\label{eq:ezittr}
r = \frac{\sqrt{1+t}}{m^2}\left( \sqrt{m^2+1-(m^2-1)t} + \sqrt{1+t} \right).
\end{equation}

Observe that $d=|r-t|$. From this and (\ref{eq:ezittr}), we have

\begin{equation}
g_m(t) = \left( \frac{\sqrt{1+t}}{m^2} \left( \sqrt{m^2+1-(m^2-1)t} + \sqrt{1+t} 
\right) - t \right)^2 +
\end{equation}
\[
+ \frac{2\sqrt{1+t}}{m^2} \left( \sqrt{m^2+1-(m^2-1)t} + \sqrt{1+t} \right) - 1.
\]

Let us observe that $g_3\left( \frac{1}{2} \right) = 0$,
and that $g_m(0) < 0$ and $g_m(1) > 0$, for every $m > 3$.
Since $g_m$ is continuous on $[0,1]$, $g_m$ has a root
in the interval $(0,1)$, for all $m \geq 3$.
Thus, for every $n \geq 4$, we have found a simplex $P$ and
a sphere $\S^{n-1}(c,r)$ that satisfy Lemma~\ref{lem:maehara counterexample} $(i)$.
\end{proof}

\section{Monotonicity of the Inradius, the Minimal Width and the Diameter of a Ball-polyhedron under a Contraction of the Centers}\label{sec:monotonicity}

One of the best known open problems of discrete geometry is the Kneser--Poulsen 
conjecture. It involves unions (resp., intersections) of finitely many balls 
in $\Re^n$ and states that, under arbitrary contraction of the center points, the 
volume of the union (resp., intersection) does not increase (resp., decrease). Recently, 
the conjecture has been proved in the plane by K. Bezdek and R. Connelly in 
\cite{BeCo} and it is open for $n\geq 3$. The interested reader is referred to the 
papers \cite{BeCo2}, \cite{Cs1}, \cite{Cs2} and \cite{Cs3} for further 
information on this problem. In this section, we investigate similar
problems.  Namely, we apply an arbitrary contraction to the center points
of the generating balls of a ball-polyhedron, and ask whether the
inradius (resp., the circumradius, the diameter and the minimum width)
can decrease.

\begin{thm}
Let $X\subset\Re^n$ be a finite point set contained in a closed unit ball of 
$\Re^n$ and let $Y$ be an arbitrary contracted image of $X$ in $\Re^n$. Then the 
inradius of $\B[Y]$ is at least as large as the inradius of $\B[X]$. 
\end{thm} 

\begin{proof}
First, observe the following fact. If $r$ denotes the inradius of $\B[X]$ (that 
is, the radius of the largest ball contained in $\B[X]$) and $R$ denotes the 
circumradius of $X$ (that is, the radius of the smallest ball containing $X$), 
then $r+R=1$. Second, we recall the following monotonicity result (see for example 
\cite{Al}). The circumradius of $X$ is at least as large as the circumradius of 
$Y$. From these two observations our theorem follows immediately.  
\end{proof}

\begin{figure}[here]
\includegraphics[width=0.5\textwidth]{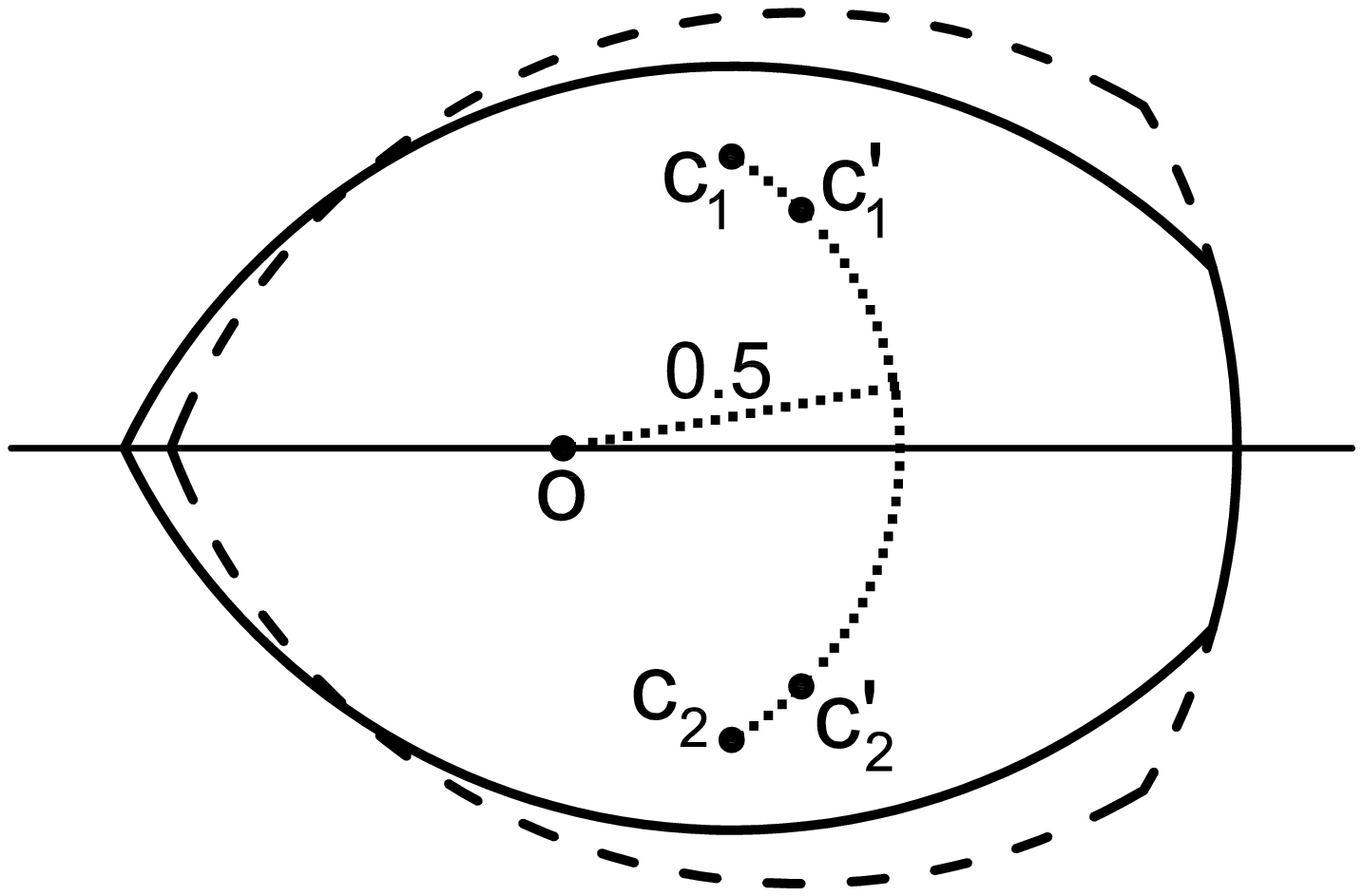}
\caption[]{}
\label{diam}
\end{figure}

The following construction (see Figure~\ref{diam}) shows that both the
diameter and the circumradius of an intersection of unit disks in the
plane can decrease under a continuous contraction of the centers. We
describe the construction in terms of polar coordinates.
The first coordinate of a vector (that is, a point) is the Euclidean
distance of the point from the origin, the second is the oriented angle
of the vector and the oriented x-axis.

Let $c_1:=(0.5,\frac{\pi}{3}), c_2:=(0.5,-\frac{\pi}{3})$,
$c_1':=(0.5,\frac{\pi}{4})$ and $c_2':=(0.5,-\frac{\pi}{4})$. Let $X$ be
the set of centers $X:=\{o,c_1,c_2\}$, and $Y:=\{o,c_1',c_2'\}$.
Clearly, $Y$ is a continuous contraction of $X$. However, a simple
computation shows that both the diameter and the circumradius of $\B[Y]$
is smaller than that of $\B[X]$.

A similar construction shows that the minimal width of an intersection
of unit disks on the plane can decrease under a continuous contraction
of the centers.

%(see Figure~\ref{diam}). However, it is easier to prove
%this using the first construction. According to the note after Theorem~1
%in \cite{Capoyleas},  

Let $c_1:=(0.8,\frac{\pi}{10}), c_2:=(0.8,-\frac{\pi}{10})$ and
$c_1':=c_2':=(0.8,0)$. Let $X$ be the set of centers
$X:=\{o,c_1,c_2\}$, and $Y:=\{o,c_1'\}$. Clearly, $Y$ is a continuous
contraction of $X$. However, a simple computation shows that the minimum
width of $\B[Y]$ is smaller than that of $\B[X]$.

%\begin{figure}[here]
%\includegraphics[width=0.3\textwidth]{minwidth.eps}
%\caption[]{}
%\label{minwidth}
%\end{figure}

\section{The Problem of Finding an Analogue to a Theorem of Steinitz for Ball-polyhedra in $\Re^3$}\label{sec:Steinitz}

One may define the vertices, edges and faces of \emph{any}
ball-polyhedron $P$ in $\Re^3$ in the natural way, as in \cite{BN}.
Henceforth in this section, we assume that $P$ is a ball-polyhedron in
$\Re^3$ with at least three balls in the reduced family of generating
balls of $P$. Then $P$ has faces that are spherically convex on the
generating sphere of $P$ that they belong to. Edges of $P$ are circular
arcs of radii less than one (not full circles) ending in vertices.
Moreover, every vertex is adjacent to at least three edges and at least
three faces of $P$.

In this paper, a \emph{graph} is non-oriented and has finitely many
vertices and edges. A graph is \emph{$2$-connected} (resp.,
\emph{$3$-connected} ), if it has at least three (resp., four) vertices
and deleting any vertex (resp., any two vertices) yields a connected graph. A graph is
\emph{simple} if it contains no loop (an edge with identical end-points)
and no parallel edges (two edges with the same two end-points).

The edge-graph of $P$ contains no loops, but may contain parallel edges.
Moreover, it is  $2$-connected and planar.

By the construction given in \cite{BN}, there is a ball-polyhedron $P$
in $\Re^3$ with two faces meeting along a series of edges. The family of
vertices, edges and faces of $P$ (together with the empty set and $P$
itself) do \emph{not} form an algebraic lattice with respect to containment.

\begin{rem}\label{rem:faceintersection}
Clearly, a ball-polyhedron $P$ in $\Re^3$ is standard if, and only if, 
the vertices, edges and faces of $P$ (together with $\emptyset$ and $P$)
form an algebraic lattice with respect to containment.

It follows that for any two faces $F_1$ and $F_2$ of a standard
ball-polyhedron $P$ in $\Re^3$, the intersection $F_1\cap F_2$ is either
empty or one vertex or one edge of $P$.
\end{rem}

In what follows, we investigate whether an analogue of the famous
theorem of Steinitz regarding the edge-graph of convex polyhedra holds for
standard ball-polyhedra in $\Re^3$. Recall that this theorem states that
a graph is the edge-graph of some convex polyhedron in $\Re^3$ if, and
only if, it is simple, planar and $3$-connected. 

\begin{claim}\label{claim:approxconvpolytopes}
Let $\bar P$ be a convex polyhedron in $\Re^3$ with the
property that every face of $\bar P$ is inscribed in a circle. Let
$\Lambda$ denote the face lattice of $\bar P$.\\
Then there is a sequence $\{P_1,P_2,\dots\}$ of standard ball-polyhedra
in $\Re^3$ with face lattices isomorphic to $\Lambda$ such that
$\lim\limits_{k\rightarrow\infty}kP_k=\bar P$ in the Hausdorff metric. 
\end{claim}

\begin{proof}
Let $\fac$ denote the set of the (two-dimensional) faces of $\bar P$;
let $c_F$ denote the circumcenter, $r_F$ the circumradius, and $n_F$ the
inner unit normal vector of the face $F\in \fac$. We define $P_k'$ as
the following intersection of closed balls of radius $k$.
\begin{equation}
P_k':=\bigcap\limits_{F\in\fac}
\B\left[c_F+\left(\sqrt{k^2-r_F^2}\right)n_F,k\right]
\end{equation}
Clearly, $P_k:=\frac{1}{k}P_k'$ is a ball-polyhedron in $\Re^3$. The
terms face, vertex and edge of $P_k'$ are defined in a natural way,
exactly as for ball-polyhedra in $\Re^3$. It is easy to see that every
vertex of $\bar P$ is a vertex of $P_k'$. Moreover, a simple
approximation argument shows that, for sufficiently large $k$, $P_k'$ is
a standard ball-polyhedron in $\Re^3$ with a face lattice that is
isomorphic to $\Lambda$. Clearly, 
$\lim\limits_{k\rightarrow\infty}P_k'=\bar P$. Now, we take a
$k_0\in\Ze^+$ such that the face lattice of $P_{k_0}$ is isomorphic to
$\Lambda$ and we replace the (finitely many) elements of the sequence
$\{P_1, P_2,\dots\}$ that have a face lattice non-isomorphic to
$\Lambda$ by $P_{k_0}$. The sequence of ball-polyhedra obtained this way
satisfies the requirements of the claim.
\end{proof}

\begin{cor}
If $\Lambda$ is a graph that can be realized as the edge-graph of a
convex polyhedron $\bar P$ in $\Re^3$, with the property that every face
of $\bar P$ is inscribed in a circle, then $\Lambda$ can be realized as
the edge-graph of a standard ball-polyhedron in $\Re^3$. 
\end{cor}

We note that \emph{not} every $3$-connected, simple, planar graph can
be realized as the edge-graph of a convex polyhedron in $\Re^3$ with all faces
having a circumcircle. See pp. 286-287 in \cite{Grun}.

\begin{claim}\label{claim:edgegraph}
The edge-graph of any standard ball-polyhedron $P$ in $\Re^3$ is
simple, planar and a $3$-connected graph.
\end{claim}

\begin{proof}
Let $G$ be the edge-graph of $P$. It is clearly planar. By a simple case
analysis, one obtains that $G$ has at least four vertices. First, we show
that $G$ is simple. Clearly, there are no loops in $G$. 

Assume that two vertices $v$ and $w$ are connected by at least two
edges $e_1$ and $e_2$.  From the reduced family of generating spheres
of $P$, let $Q$ be the intersection of those that contain $e_1$ or
$e_2$. Clearly, $Q=\{v,w\}$ which contradicts
Remark~\ref{rem:faceintersection}.

Now, we show that $G$ is $3$-connected. Let $v$ and $w$ be two arbitrary
vertices of $G$.  Take two vertices $s$ and $t$ of $G$, both different
from $v$ and $w$. We need to show that there is a path between $s$ and
$t$ that avoids $v$ and $w$. We define two subgraphs of $G$, $C_v$ and
$C_w$ as follows. Let $C_v$ (resp., $C_w$) be the set of vertices of $P$
that lie on the same face as $v$ (resp., $w$) and are distinct
from $v$ (resp., $w$). Let an edge $e$ of $G$ connecting two points of
$C_v$ (resp., $C_w$) be an edge of $C_v$ (resp., $C_w$) if, and only if,
$e$ is an edge of a face that contains $C_v$ (resp., $C_w$).

By Remark~\ref{rem:faceintersection}, $C_v$ and $C_w$ are cycles.
Moreover, $v$ and $w$ are incident to at most two faces in common.

Case 1: $v$ and $w$ are not incident to any common face; that is,
$v\notin C_w$ and $w\notin C_v$. Since $G$ is connected, there is a path
connecting $s$ and $t$. We may assume that this path does not pass
through any vertex twice. Assume that this path includes $v$ by passing
through two edges, say $e_1$ and $e_2$ that share $v$ as a vertex. Let
the vertex of $e_1$ (resp., $e_2$) different from $v$ be $v_1$ (resp.,
$v_2$). Clearly, $v_1,v_2\neq w$ and they are contained in $C_v$, which
is a cycle. Thus, the edges $e_1$ and $e_2$ in the path may be replaced
by a sequence of edges of $C_v$ that connects $v_1$ and $v_2$. If the
path passes through $w$ then it may be modified in the same manner to
avoid $w$, thus we obtain the desired path. 

Case 2: $v$ and $w$ are incident to one or two common faces. Let
$C$ be the subgraph of $C_v\cup C_w$ spanned by the union of vertices of
$C_v$ and $C_w$ erasing $v$ and $w$. Since $P$ is a standard
ball-polyhedron, $C$ is a cycle. Similarly to the preceding argument,
any path from $s$ to $t$ may be modified such that it does not pass
through $v$ and $w$ using edges of $C$.
\end{proof}

We pose the following questions.

\begin{prob}
Prove or disprove that every $2$-connected planar graph with no loops is the
edge-graph of a ball-polyhedron in $\Re^3$.
\end{prob}

\begin{prob}
Prove or disprove that every $3$-connected, simple, planar graph is
the edge-graph of a standard ball-polyhedron in $\Re^3$.
\end{prob}

\section{Ball-polyhedra in $\Re^3$ with Symmetric Sections}\label{sec:symmetricsections}

Let $K\subset\Re^3$ be a convex body with the property that any planar
section of $K$ is axially symmetric. The first named author conjectured
(see \cite{GrOd}) that, in this case, $K$ is either a body of revolution
or an ellipsoid. A remarkable result related to the conjecture is due to
Montejano \cite{Montejano}. He showed that if $K\subset\Re^3$ is a
convex body with the property that, for some point $p\in\inter K$, every
planar section of $K$ through $p$ is axially symmetric, then there is a
planar section of $K$ through $p$ which is a disk. Unfortunately, the
claim of \'Odor  (see \cite{GrOd}) that he proved this conjecture turned
out to be too optimistic, his approach was found incomplete.  The
following theorem shows that the conjecture holds for the class of
ball-polyhedra in $\Re^3$ with the weaker condition in Montejano's
result.

\begin{thm}
Let $P$ be a ball-polyhedron in $\Re^3$ and $p\in\inter P$ with the property
that any planar section of $P$ through $p$ is axially symmetric.\\
Then $P$ is either one point or a unit ball
or the intersection of two unit balls.
\end{thm}
\begin{proof}
Assume the contrary; that is, that the minimum number of unit balls
needed to generate $P$ is $k\geq 3$. Let the reduced family of
generating unit spheres be $\S^{n-1}(c_1),\dots,\S^{n-1}(c_k)$. Since
$P$ is generated by at least three unit balls, it has an edge. Let $q_1$
be any point in the relative interior of some edge $e$ of $P$ and let
$q_2$ be a point in the relative interior of a facet $F$ of $P$ that
does \emph{not} contain $e$.  By slightly moving $q_1$ on $e$ and $q_2$
on $F$, we may assume that the plane $H$ spanned by $p, q_1$ and $q_2$
does not contain any vertex of $P$ and is neither parallel nor
perpendicular to the line passing through $c_i$ and $c_j$, for any
$1\leq i< j\leq k$.

Since $F$ does not contain $e$, it follows that $H$ intersects at least
three edges of $P$. Thus, $H\cap P$ is a convex planar figure in $H$
bounded by a closed curve that is a series of at least three circular
arcs. Moreover, since $H$ is neither parallel nor perpendicular to the
line passing through $c_i$ and $c_j$, for any $1\leq i< j\leq k$, the
radii of these arcs are pairwise distinct. This clearly contradicts our
assumptions on $P$ as such a planar figure is \emph{not} axially symmetric.
\end{proof}
                                                                                                                  
\section{Illumination of Ball-polyhedra and Sets of Constant Width in $\Re^3$}\label{sec:illumination}

We consider the illumination problem for ball-polyhedra in $\Re^3$ that
contain the centers of their generating balls. We prove that such bodies
are illuminated by three pairs of opposite directions that are mutually
orthogonal. The method we use naturally extends to bodies obtained as
intersections of infinitely many balls, hence it yields a proof of the
known theorem (Lassak \cite{Lassak}, and Weissbach
\cite{Weissbach}) that any set of constant width in $\Re^3$ is
illuminated by six light sources. For a survey on illumination see
\cite{MartiniSoltan} and the new paper of K. Bezdek \cite{B}; bodies of constant 
width are discussed in the surveys 
\cite{ChakGro} and \cite{MSW}, see also the monograph \cite{YagBol}.

\begin{defn}
 Let $K\in\Re^n$ be a convex body and $z\in\bd K$ a point on its
 boundary. We say that the direction $\vu\in\Sph$ \emph{illuminates} $K$
 at $z$ if the ray $\{z+t\vu : t>0\}$ intersects the interior of $K$. 
 Furthermore, $K$ is \emph{illuminated} at $z\in\bd K$ \emph{by a set
 $A\subseteq\Sph$ of directions} if at least one direction from $A$
 illuminates $K$ at $x$. Then $K$ is illuminated by $A\subset\Sph$ if
 $K$ is illuminated by $A$ at every boundary point of $K$.
\end{defn}

Let $K$ be a convex body and $z\in\bd K$. We denote by $G(z)$ the
set of inward unit normal vectors of hyperplanes that support $K$ at
$z$. We note that $-G(z)$ is the Gauss image of $z$.

We denote the open hemisphere of $\Sph$ with center $\vu\in\Sph$ by
$D(\vu)$ and its relative boundary (a great sphere of $\Sph$) by
$C(\vu)$. Then $\vu$, which is in $\Sph$, illuminates $K$ at $z\in\bd
K$ if, and only if, $G(z)\subset D(\vu)$. This leads to the following
observation which is an easy special case of the Separation Lemma in
\cite{BIllumin}.

\begin{obs}\label{obs:blocking}
The pair of directions $\{\pm\vu\}\subset\Sph$ illuminates the convex
body $K\subset\Re^n$ at $z\in\bd K$ if, and only if, $G(z)\cap
C(\vu)=\emptyset$.
\end{obs}

Now, we are ready to state the main result of this section.
\begin{thm}\label{thm:illumin}
Let $X\subset\Re^3$ be a set of diameter at most one, and let
$\vu\in\S^2(o)$ be given.\\  Then there exist $v$ and $w$ in $\Sph$
such that $u, v$ and $w$ are pairwise orthogonal, and the body $K:=\B[X]$
is illuminated by the six directions $\{\pm u,\pm v, \pm w\}$.
\end{thm}

Using the fact that a closed set $X\subset\Re^n$ is of constant width
one if, and only if, $\B[X]=X$, cf. \cite{Egg}, we obtain the following corollary.

\begin{cor}
Any set of constant width can be illuminated by three pairwise orthogonal pairs
of opposite directions, one of which can be chosen arbitrarily.
\end{cor}

To prove the theorem we need the following lemma.
\begin{lem}\label{lem:Gaussimage}
Let $X\subset\Re^n$ be a set of diameter at most one and $z\in\bd \B[X]$.\\
Then $G(z)\subset\Sph$ is of spherical diameter not greater than
$\frac{\pi}{3}$.
\end{lem}
\begin{proof}
We may assume that $X$ is closed. It is easy to see that\\
$G(z)=\Sconv{\S^{n-1}(z)\cap X}{\S^{n-1}(z)}-z$. So, we have to
show that if $x_1, x_2\in\S^{n-1}(z)\cap X$, then
$\sphericalangle(x_1zx_2)\leq\frac{\pi}{3}$. It is true, since the
Euclidean isosceles triangle $\conv\{x_1,z,x_2\}$ has two legs $[z,x_1]$ and $[z,x_2]$
of length one, and base $[x_1,x_2]$ of length at most one,
because the diameter of $X$ is at most one. This proves the lemma.
\end{proof}

\begin{proof}[Proof of Theorem \ref{thm:illumin}]
Let the direction $\vu\in\S^2(o)$ be given. We will call $\vu$ vertical,
and directions perpendicular to $\vu$ horizontal. We pick two pairwise
orthogonal, horizontal directions, $\vv_1$ and $\vv_2$. Assume that the
six directions $\{\pm \vu, \pm \vv_1, \pm \vv_2\}\subset\S^2$ do not
illuminate $K$. According to Observation \ref{obs:blocking}, there is a point
$z\in\bd K$ such that $G(z)$ intersects each of the three great circles
of $\S^2(o)$: $C(\vu), C(\vv_1)$ and $C(\vv_2)$. We choose three points
of $G(z)$, one on each great circle: $y_0\in G(z)\cap C(\vu),y_1\in
G(z)\cap C(\vv_1)$ and $y_2\in G(z)\cap C(\vv_2)$. Note that each of the
three great circles is dissected into four equal arcs (of length
$\frac{\pi}{4}$) by the two other great circles.

By Lemma \ref{lem:Gaussimage}, $G(z)\subset\S^2(o)$ is a spherically
convex set of spherical diameter at most $\frac{\pi}{3}$. However,
$y_0,y_1,y_2\in G(x)$, so the generalization of Jung's theorem for
spherical space by Dekster \cite{Dekster} shows that $y_0,y_1$ and $y_2$
are the mid-points of the great circular arcs mentioned above. So, the
only way that a point $z\in\bd K$ is not illuminated by any of the six
directions $\{\pm \vu, \pm \vv_1, \pm \vv_2\}$ is the following. The set
$G(z)$ contains a spherical equilateral triangle of spherical side
length $\frac{\pi}{3}$ and the vertices of this spherical triangle lie
on $C(u), C(v_1)$ and $C(v_2)$, respectively. Furthermore, each vertex is
necessarily the mid-point of the quarter arc of the great circle on
which it lies, and $G(z)$ does not intersect either of the three great
circles in any other point.

Since the set $\{G(z) : z\in\bd K\}$ is a tiling of $\S^2(o)$, there are only
finitely many boundary points $z\in\bd K$ such that $G(z)$ contains an
equilateral triangle of side length $\frac{\pi}{3}$ that has a vertex
on $C(\vu)$. We call these tiles blocking tiles.

Now, by rotating $\vv_1$ and $\vv_2$ together in the horizontal plane,
we can easily avoid all the blocking tiles; that is, we can find a rotation $R$
about the line spanned by $\vu$ such that none of the blocking tiles has a vertex
on both circles $C(R(\vv_1))$ and $C(R(\vv_2))$.
Now, $\pm \vu, \pm R(\vv_1)$ and $ \pm R(\vv_2)$ are the desired
directions finishing the proof of the theorem.
\end{proof}

We remark that in the theorem, we can ``almost'' choose the second
direction arbitrarily, more precisely: Given any two orthogonal vectors
$\vu, \vv_1\in\S^2(o)$ and $\varepsilon>0$, we may find two directions
$\vv_1', \vv_2'\in\S^2(o)$ such that $\pm \vu,\pm \vv_1'$ and $\pm
\vv_2'$ illuminate $K$, $vu, \vv_1$ and $\vv_2$ are pairwise orthogonal,
and $\|\vv_1-\vv_1'\|<\varepsilon$.\\
This statement may be derived from the last paragraph of the proof. The
set of rotations about the line spanned by $\vu$ is a one-parameter
group parametrized by angle. Each blocking tile rules out at most four
angles, and there are finitely many blocking tiles. This argument proves
the following statement.

\begin{thm}
Let $X\subset\Re^3$ be a set of diameter at most one. We choose three
pairwise orthogonal directions $u, v$ and $w$ in $\S^2(o)$ randomly with
a uniform distribution.\\
Then the body $K:=\B[X]$ is illuminated by $\{\pm u, \pm v, \pm w\}$
with probability one.
\end{thm}

\begin{prob}
Let $X\subset\Re^3$ be a set of diameter at most one. Prove or disprove
that $\B[X]$ is illuminated by four directions.
\end{prob}

\section{Dowker-Type Isoperimetric Inequalities for Disk-Polygons}\label{sec:Dowker}

In this section, we examine theorems concerning disk-polygons that are
analogous to those studied by Dowker in \cite{Do} and L. Fejes T\'oth in
\cite{FT} for polygons. The arguments are based on (\cite{FT}, pp.162-170), but are
adapted to the current setting using \cite{BeCoCs} and \cite{CsLN}.

Let $x_0, x_1, \ldots, x_n$ be points in the plane such that
they are all distinct, except for $x_0$ and $x_n$, which are equal.
Furthermore, suppose that the distance between each pair of consecutive
points is at most two. Next, let $\widehat{x_ix_{i+1}}$ denote one of
the two unit circle arcs of length at most $\pi$ with endpoints $x_i$
and $x_{i+1}$. A \emph{circle-polygon} is the union of these unit circle
arcs, $\widehat{x_0x_1}, \widehat{x_1x_2}, \ldots,
\widehat{x_{n-1}x_n}$. The points $x_0, x_1, \ldots, x_n$ are the
\emph{vertices} of the circle-polygon and the unit circle arcs
$\widehat{x_0x_1}, \widehat{x_1x_2}, \ldots, \widehat{x_{n-1}x_n}$ are
the \emph{edges}, or more commonly the \emph{sides}, of the
circle-polygon.  Finally, the \emph{underlying polygon} is the polygon
formed by joining the vertices, in order, by straight line segments.
Observe that both a circle-polygon and its underlying polygon may have
self-intersections.

The definition of a standard ball-polytope implies that a disk-polygon
is standard if, and only if, there are at least three disks in the
reduced family of generating disks. Such a disk-polygon $P$ has well
defined vertices and edges. Clearly, $\bd P$ is a circle-polygon with
vertices and edges which coincide with those of $P$. The underlying
polygon of an $n$-sided  disk-polygon is just the boundary of the convex
hull of the vertices. An $n$-sided  disk-polygon is called
\emph{regular} if the underlying polygon is regular.

Let $C$ be a circle of radius $r < 1$. A circle-polygon (resp.
disk-polygon) $P$ is \emph{inscribed in} $C$ if $P \subset \conv C$
and the vertices of $P$ lie on $C$. A circle-polygon (resp.
disk-polygon) $P$ is \emph{circumscribed about} $C$ if $C \subset P$
and the interior of each edge of $P$ is tangent to $C$. A standard
compactness argument ensures the existence of an $n$-sided disk-polygon
of largest (resp. smallest) perimeter, as well as one of largest (resp.
smallest) area, inscribed in (resp. circumscribed about) $C$.

\begin{lem}\label{lem:concper}
Let $C$ be a circle of radius $r < 1$. Let $P_n$ be an $n$-sided disk-polygon of largest perimeter inscribed in $C$. Then
\begin{equation}
\per(P_{n-1}) + \per(P_{n+1}) < 2 \per(P_n), \mbox{ for all } n\geq4.
\end{equation}
\end{lem}

\begin{proof}
Let $Q$ be an $(n-1)$-sided  disk-polygon and $R$ be an $(n+1)$-sided disk-polygon, both inscribed in $C$. To prove the theorem we need only construct two $n$-sided  disk-polygons $S$ and $T$ such that

\begin{equation}\label{eqn:specialconcper}
\per(Q) + \per(R) \leq \per(S) + \per(T).
\end{equation}

Without loss of generality we make the following assumptions. First, inscribe $Q$ and $R$ into $C$ so that their respective vertices do not coincide. Second, any arc of $C$ with length at least $\pi r$ contains a vertex from each of $Q$ and $R$. Otherwise, there exists an $(n-1)$-sided disk-polygon (resp. $(n+1)$-sided  disk-polygon) with larger perimeter than $Q$ (resp. $R$).

Let $x_i$ and $x_{i+1}$ be two consecutive vertices of a circle-polygon
$P$ inscribed in $C$. Suppose that one of the arcs of $C$ from $x_i$ to
$x_{i+1}$ contains neither $x_{i-1}$ nor $x_{i+2}$. Let $\bar{C}$ denote
this arc. A \emph{cap} of $C$ from $x_i$ to $x_{i+1}$, denoted by
$C(x_i,x_{i+1})$, is the segment of $\conv C$ bounded by $\bar{C}$ and
the line segment through $x_i$ and $x_{i+1}$. Now, suppose that there is
a cap contained in another cap (see Figure~\ref{fig:dow1}). More
precisely, there are four vertices $a, b, l, m$ such that the vertices
$a$ and $b$ (resp. $l$ and $m$) form an edge $\widehat{ab}$ (resp.
$\widehat{lm}$) and $C(l,m) \subset C(a,b)$.

\begin{figure}[here]
\includegraphics[width=0.5\textwidth]{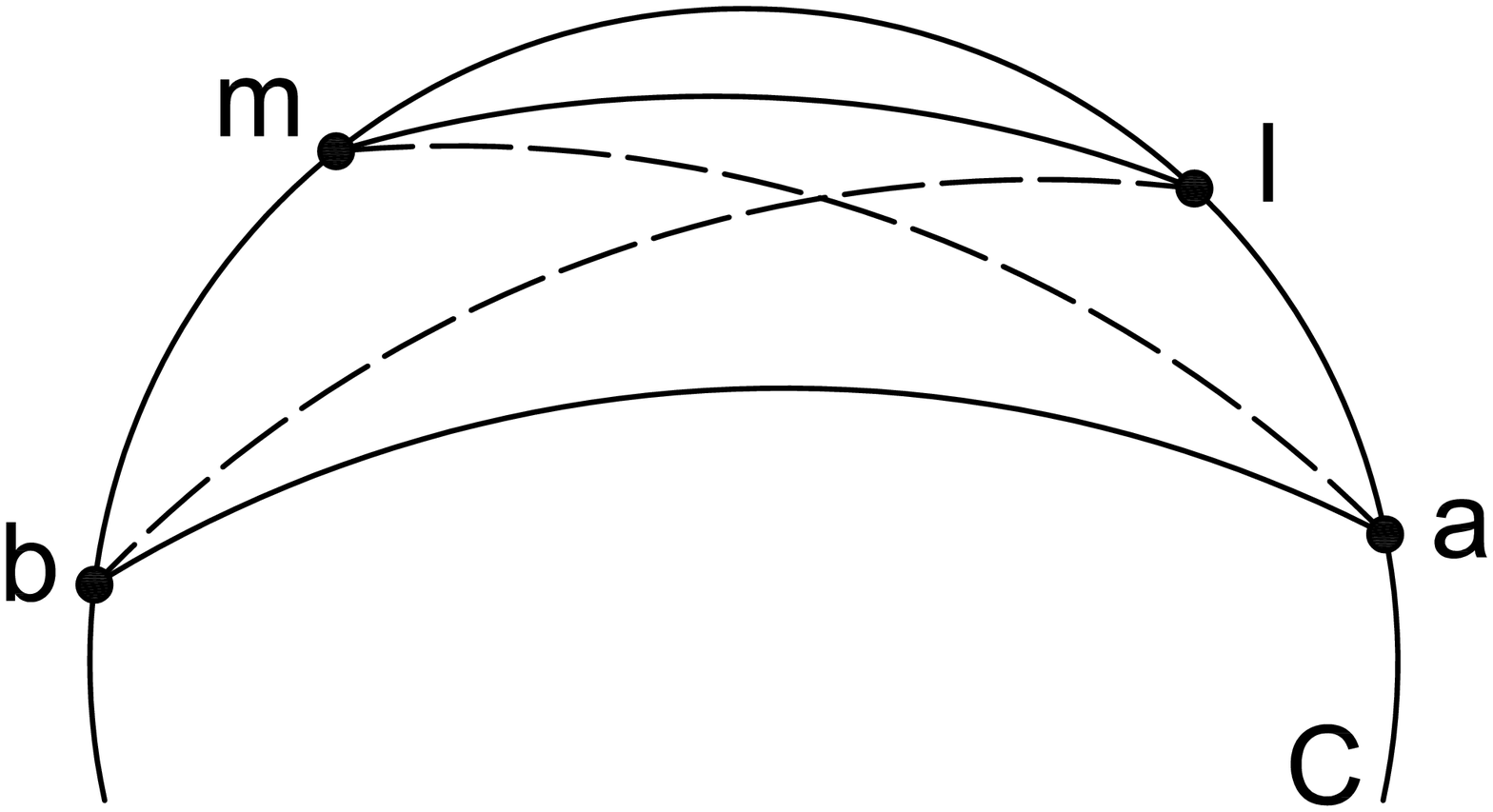}
\caption[]{}
\label{fig:dow1}
\end{figure}

This configuration, where a cap is contained in another cap, may arise
when two circle-polygons, say $A$ and $B$, are inscribed in a single
circle, or when a single self-intersecting circle-polygon, say $U$, is
inscribed in a circle. Suppose that we start with the former (resp.
latter) and as in the figure, assume that the cyclic ordering of the
vertices is $a, l, m, b$. Choosing either of the unit circle arcs
joining $a$ to $m$, we obtain $\widehat{am}$. Similarly, we construct
$\widehat{bl}$. After replacing the edges $\widehat{ab}$ and
$\widehat{lm}$ with $\widehat{am}$ and $\widehat{bl}$, respectively, we
obtain a single self-intersecting circle-polygon, which we call $U$
(resp. two circle-polygons, which we call $A$ and $B$). By the
inequality mentioned in Corollary \ref{cor:circlequadrilateral}, the
total perimeter of $U$ (resp. $A$ and $B$) is strictly larger than the
total perimeter of $A$ and $B$ (resp. $U$).

Starting from the circle-polygons $Q$ and $R$, we carry out the
preceding algorithm for each cap contained in another cap. After every
odd numbered iteration of the algorithm we obtain a single,
self-intersecting circle-polygon and after every even numbered iteration
we obtain two circle-polygons with no self-intersections. Since we have
finitely many vertices and the perimeter increases strictly with each
step, the algorithm terminates. Furthermore, when it does terminate,
there is no cap contained within another cap.

A simple counting argument shows that the process terminates with two circle-polygons. Each one is an $n$-sided  circle-polygon because no cap of one is contained within a cap of the other. Let us denote these $n$-sided  circle-polygons by $S^*$ and $T^*$. Since the perimeter increased at each step of the process,

\begin{equation}
\per(Q) + \per(R) < \per(S^*) + \per(T^*).
\end{equation}

Now, it may be the case that some of the edges of $S^*$ (resp. $T^*$) have relative interior points that meet the convex hull of the underlying polygon. We replace each such edge by the other shorter unit circle arc passing through the same two vertices. It is clear that this produces a disk-polygon with the same perimeter as $S^*$ (resp. $T^*$). This is the desired $S$ (resp. $T$).
\end{proof}

\begin{thm}
Let $C$ be a circle of radius $r<1$. Let $P$ be an $n$-sided disk-polygon of largest perimeter that can be inscribed in $C$.
Then $P$ is regular.
\end{thm}

\begin{proof}
Suppose that $P$ is not regular. Starting with $P$ and a suitable rotation of $P$ we modify the argument in the proof of the preceding lemma to construct two $n$-sided  disk-polygons $Q$ and $R$ inscribed in $C$. By construction, $\per(Q) + \per(R) > 2 \per(P)$. Hence, one of $Q$ or $R$ has larger perimeter than $P$.
\end{proof}

\begin{thm}\label{thm:in_max_area}
Let $C$ be a circle of radius $r<1$. Let $P$ be an $n$-sided disk-polygon of largest area that can be inscribed in $C$. Then $P$ is
regular.
\end{thm}

\begin{proof}
Suppose $P$ is not regular. Let $P_0$ be the regular $n$-sided
disk-polygon with the same perimeter as $P$. By the discrete
isoperimetric inequality for circle-polygons proved in \cite{CsLN},
$\area(P) < \area(P_0)$. Furthermore, by the preceding theorem, $P_0$ is
inscribed in a circle $C_0$ with radius $r_0 < r$. Thus, $P_1$, the
regular $n$-sided disk-polygon inscribed in $C$, clearly satisfies
$\area(P_0) < \area(P_1)$ which completes the proof.
\end{proof}

In general, the behavior of the areas of disk-polygons inscribed in a circle is difficult to describe, but we do so in the following special case.

\begin{lem}
Let $C$ be a circle of radius $r < 1$. Let $P_n$ be an $n$-sided disk-polygon of largest area inscribed in $C$. Then
\begin{equation}
\area(P_{n-1}) + \area(P_{n+1}) < 2 \area(P_n), \mbox{ for all odd $n$ and } n\geq5.
\end{equation}
\end{lem}

\begin{proof}
By Theorem~\ref{thm:in_max_area}, $P_{n-1}$ and $P_{n+1}$ are regular. Since $n$ is odd, both $P_{n-1}$ and $P_{n+1}$ are symmetric about a line through opposite vertices. After an appropriate rotation, $P_{n-1}$ and $P_{n+1}$ share such a line of symmetry $d$ which separates the symmetric sections of each of $P_{n-1}$ and $P_{n+1}$ (see Figure~\ref{fig:dow2}).

\begin{figure}[here]
\includegraphics[width=0.6\textwidth]{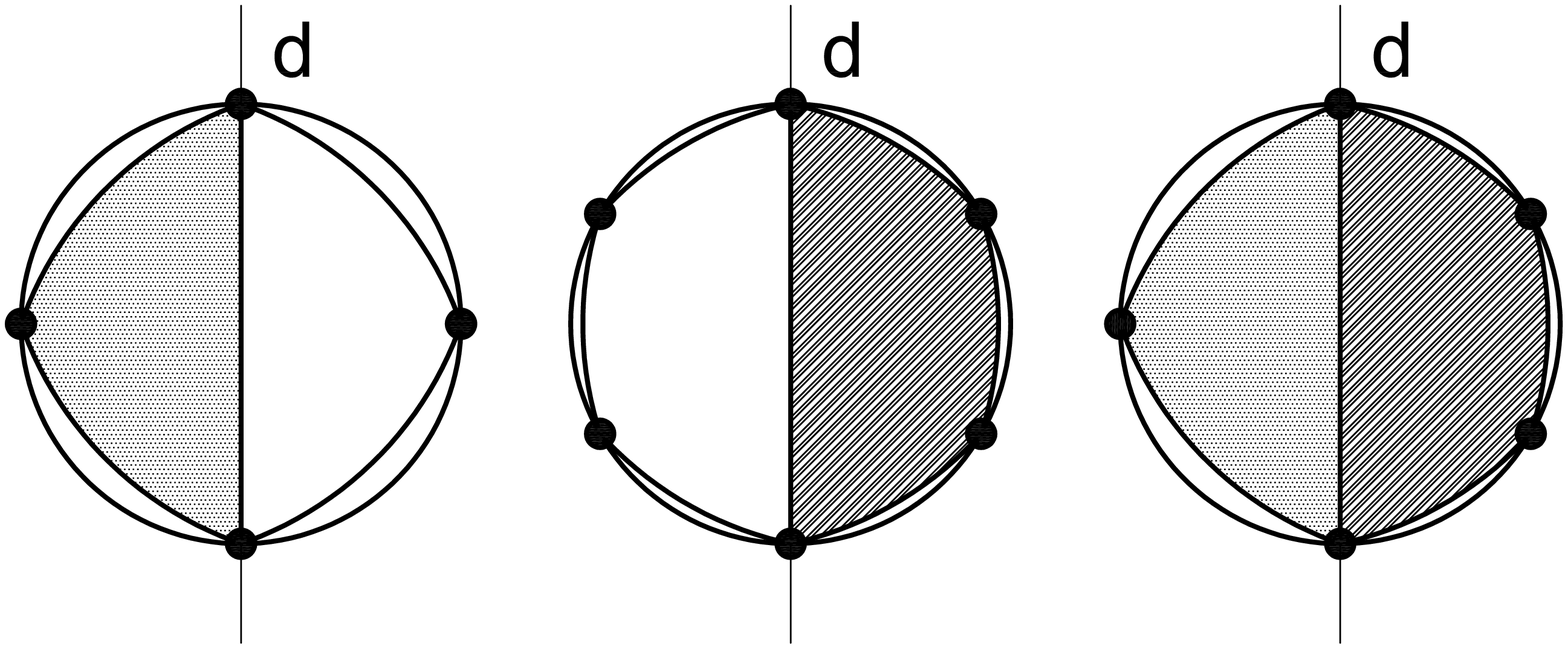}
\caption[]{}
\label{fig:dow2}
\end{figure}

Let $Q$ be one half of $P_{n-1}$ lying on one side of $d$ and $R$ the half of $P_{n+1}$ lying on the other side of $d$. The union of $Q$ and $R$ is an $n$-sided  disk-polygon $U$ inscribed in $C$. Clearly, $\area(P_{n-1}) + \area(P_{n+1}) = 2 \area{U} < 2 \area(P_n)$.
\end{proof}

There is no straight-forward method to generalize this result to all $n$-sided disk-polygon. The method described does not apply and no formula
is known to describe this area. So, we make the following conjecture.

\begin{conj}
Let $C$ be a circle of radius $r < 1$. Let $P_n$ be an $n$-sided disk-polygon of largest area inscribed in $C$. Then
\begin{equation}
\area(P_{n-1}) + \area(P_{n+1}) < 2 \area(P_n), \mbox{ for all } n\geq4.
\end{equation}
\end{conj}

We now turn our attention to disk-polygons which are circumscribed about a circle. A modification of \cite{FT}, p. 163, similar to the one described above provides the following theorem.

\begin{thm}
Let $C$ be a circle of radius $r < 1$.

\noindent (i) Let $P_n$ be an $n$-sided  disk-polygon of smallest area circumscribed about $C$. Then
\begin{equation}\label{eqn:out_area}
\area(P_{n-1}) + \area(P_{n+1}) > 2 \area(P_n), \mbox{ for all } n\geq4.
\end{equation}
\noindent Furthermore, $P_i$ is regular for all $i\geq3$.

\noindent (ii) Let $P_n$ be an $n$-sided  disk-polygon of smallest perimeter circumscribed about $C$. Then
\begin{equation}\label{eqn:out_per}
\per(P_{n-1}) + \per(P_{n+1}) > 2 \per(P_n), \mbox{ for all } n\geq4.
\end{equation}
\noindent Furthermore, $P_i$ is regular for all $i\geq3$.
\end{thm}

\section{Erd\H os--Szekeres type Problems for Ball-polytopes}\label{sec:Erdos-Szekeres}

\begin{defn}
Let $A \subset \Re^n$ be a finite set contained in a closed unit ball.
Then $\convv A$ is called a ball-polytope.
\end{defn}

In this section we are going to find analogues of some results about
convex polytopes for ball-polytopes.
We begin with two definitions.

\begin{defn}
Let $A \subset \Re^n$ be a finite set.
If $x \notin \conv(A \setminus \{x\})$ for any $x \in A$,
we say that the points of $A$ are in convex position.
\end{defn}

\begin{defn}\label{defn:EuclESZ}
For any $n \geq 2$ and $m \geq n+1$, let $f_n(m)$ denote the maximal cardinality
of a set $A \subset \Re^n$ that satisfies the following two conditions:

\begin{tabular}{cl}
(i)  & any $n+1$ points of $A$ are in convex position,\\
(ii) & $A$ does not contain $m$ points that are in convex position.
\end{tabular}
\end{defn}

In \cite{ESZ1} and \cite{ESZ2}, Erd\H os and Szekeres
proved the existence of $f(m) := f_2(m)$ for every $m$,
and gave the estimates $2^{m-2} \leq f(m) \leq {2m-4 \choose m-2}$.
They conjectured that $f(m) = 2^{m-2}$.
Presently, the best known upper bound is $f(m) \leq {2m-5 \choose m-2} + 1$,
given by T\'oth and Valtr in \cite{TV}.
We note that, if the projections of $m$ points of $\Re^n$ to an affine subspace
are in convex position, then the original points are also in convex position.
Thus, the results about $f_2(m)$ imply also that $f_n(m)$ exists, and
$f_{n+1}(m) \leq f_n(m)$, for every $n$ and $m$.

\begin{defn}
Let $A \subset \Re^n$ be a finite set contained in a closed unit ball.
If $x \not\in \convv (A \setminus \{x\})$, for every $x \in A$, then we say that the
points of $A$ are in spindle convex position.
\end{defn} 

\begin{defn}\label{defn:spindle convESZ}
For $n \geq 2$ and $m \geq n+1$, let $g_n(m)$ be the maximal cardinality of a set
$A \subset \Re^n$ that is contained in a closed unit ball and satisfies the following properties:

\begin{tabular}{cl}
(i)  & any $n+1$ points of $A$ are in spindle convex position,\\
(ii) & $A$ does not contain $m$ points in spindle convex position.
\end{tabular}
\end{defn}

To show the importance of $(i)$ in Definition~\ref{defn:spindle convESZ},
we provide the following example.
Let $A := \{ x_1, x_2, \ldots, x_k \}$, where  $x_1, x_2, \ldots, x_k$ are points
of an arc of radius $r>1$ in this cyclic order.
Then any $n+1$ points of $A$ are affine independent whereas $A$ does not contain three points
in spindle convex position.

In the remaining part of this section, we show that $f_n(m) = g_n(m)$, for every $n$ and $m$.
Let us assume that $A \subset \Re^n$ is a set that satisfies $(i)$ and $(ii)$ in
Definition~\ref{defn:EuclESZ}.
Observe that, for a suitably small $\varepsilon>0$, any $n+1$ points of
$\varepsilon A$ are in spindle convex position.
This implies that $f_n(m) \leq g_n(m)$.
To show the inequality $f_n(m) \geq g_n(m)$, we prove the following
stronger version of Theorem~\ref{thm:CaratheodorySteinitz}.

\begin{thm}
Let $P \subset \Re^n$ be contained in a closed unit ball,
and $p \in \convv P$.
Then $p \in \conv P$ or $p \in \convv Q$, for some $Q \subset P$ with $\card Q \leq n$.
\end{thm}

\begin{proof}
We show that if $p \notin \convv Q$, for any $Q \subset P$ with $\card Q \leq n$,
then $p \in \B^n[c,r]$, for any ball $\B^n[c,r]$ that contains $P$.
Since $\conv P$ is the intersection of all the balls that contain $P$ and have radii at least one,
this will imply our statement.
We assume that there is a ball $B^n[q,r]$ with $r \geq 1$ that
contains $P$ but does not contain $p$.

If $p \in \bd \convv P$, our statement follows from Theorem~\ref{thm:CaratheodorySteinitz}.
So, let us assume that $p \in \inter \convv P$.
From this and Lemma~\ref{lem:supportbyball}, we have $p \in \inter \B^n[c,1] = \B^n(c,1)$
whenever $P \subset \B^n[c,1]$.
If, for every $r > 1$, there is a ball $\B^n[c_r,r]$ that contains $P$ but does not contain $p$,
then Blaschke's Selection Theorem guarantees the existence of a unit ball $\B^n[c,1]$
such that $P \subset \B^n[c,1]$ and $p \notin \B^n(c,1)$, a contradiction.
So, there is an $r > 1$ such that $P \subset \B^n[c,r]$ implies that $p \in \B^n[c,r]$.
Clearly, if $1 < r_1 < r_2$ and $r_2$ satisfies this property then $r_1$ also satisfies it.
Thus, there is a maximal value $R$ satisfying this property. 
Corollary~\ref{cor:intersection} suggests the notation

\begin{equation}
P(r) := \bigcap \{ \B^n(c,r) : P \subset \B^n(c,r)\}.
\end{equation}

Observe that $P(r_2) \subset P(r_1)$, for every $1 < r_1 < r_2$, and that $p \in \bd P(R)$.
Hence, applying Corollary~\ref{cor:intersection} and Theorem~\ref{thm:CaratheodorySteinitz}
for $\frac{1}{R} P$, we obtain a set $Q \subset P$ of cardinality at most $n$
such that any ball of radius $R$ that contains $Q$ contains also $p$.
We define $Q(r)$ similarly to $P(r)$ and so we have $Q(R) \subset Q(1) = \convv Q$,
which implies our statement.
\end{proof}

So, if $P \subset \Re^n$ is contained in a closed unit ball,
$\card P > f_n(m)$ and any $n+1$ points of $P$ are in spindle convex position,
then $P$ contains $m$ points in convex position, which, according to our theorem,
are in spindle convex position.

We note that our theorem implies the spindle convex analogues
of numerous other Erd\H os--Szekeres type results.
As examples, we mention \cite{AKP}, \cite{BH} and \cite{V}.

\noindent
{\bf Acknowledgement.}
The authors thank both referees their valuable comments that helped
to improve the level of presentation of the results in this paper.

\end{document}